\documentclass[11pt]{amsart}

\usepackage{amsmath,amssymb,amsthm,mathtools}
\usepackage{tikz,tikz-cd,placeins}
\usepackage[margin=1in]{geometry}
\usepackage[colorlinks=true,linkcolor=blue,citecolor=blue,urlcolor=blue]{hyperref}
\usetikzlibrary{calc,shapes.geometric}

\definecolor{d5blue}{HTML}{315E9A}
\definecolor{ns5red}{HTML}{B63D4B}
\definecolor{ink}{HTML}{17202A}

\newtheorem{theorem}{Theorem}[section]
\newtheorem{proposition}[theorem]{Proposition}
\newtheorem{lemma}[theorem]{Lemma}
\newtheorem{corollary}[theorem]{Corollary}
\newtheorem{definition}[theorem]{Definition}
\newtheorem{remark}[theorem]{Remark}
\newtheorem{example}[theorem]{Example}

\newtheorem{question}[theorem]{Question}

\newtheorem{mainthm}{Theorem}

\newcommand{\RR}{\mathsf R}
\newcommand{\KK}{\mathsf K}
\newcommand{\R}{\mathbf R}
\newcommand{\C}{\mathbf C}
\newcommand{\Q}{\mathbf Q}
\newcommand{\Z}{\mathbf Z}
\newcommand{\PP}{\mathbf P}
\newcommand{\T}{\mathbf T}
\newcommand{\A}{\mathbf A}

\newcommand{\X}{\mathfrak X}
\newcommand{\Ah}{\widehat A}
\newcommand{\Hom}{\operatorname{Hom}}
\newcommand{\End}{\operatorname{End}}
\newcommand{\im}{\operatorname{im}}
\newcommand{\lch}{\operatorname{lch}}

\newcommand{\dual}{!}
\newcommand{\h}{\hbar}
\newcommand{\bs}{{\color{blue}\char`\\}}
\newcommand{\fs}{{\color{red}/}}
\newcommand{\ttt}[1]{\text{\ttfamily #1}}
\newcommand{\Stab}{\operatorname{Stab}}
\newcommand{\Att}{\operatorname{Att}}
\newcommand{\wt}{\operatorname{wt}}
\newcommand{\Mix}{\operatorname{Mix}}
\newcommand{\gr}{\operatorname{graph}}
\newcommand{\ann}{\operatorname{ann}}
\newcommand{\ol}[1]{\overline{#1}}
\newcommand{\pair}[2]{\langle #1,#2\rangle}
\newcommand{\res}[1]{\left.#1\right|}
\newcommand{\cW}{\mathcal W}
\newcommand{\cG}{\mathcal G}

\title[GM duality for affine type $A$ bow varieties]
{Goresky--MacPherson duality \\ for bow varieties of affine type $A$}

\author{Junzhe Lyu}
\address{Department of Mathematics, University of North Carolina at Chapel Hill, USA}
\email{taiat@unc.edu}

\author{Rich\'ard Rim\'anyi}
\address{Department of Mathematics, University of North Carolina at Chapel Hill, USA}
\email{rimanyi@email.unc.edu}

\begin{document}

\begin{abstract}
Goresky and MacPherson observed that for certain pairs of holomorphic
symplectic varieties, the arrangements of linear subspaces cut out by restriction maps to torus fixed points in degree-two equivariant cohomology are dual to each
other. This \emph{GM duality} is among the expected features of a 3d mirror
(symplectic dual) pair.  We prove GM duality for all 3d mirror pairs of bow
varieties of affine type $A$, in the form natural to 3d mirror symmetry.
Namely, we prove it on the tautological model of degree-two equivariant
cohomology, spanned by the equivariant parameters and the determinants of
the tautological bundles.
We then show a consequence for the poset structures induced by attracting sets on 3d mirror pairs.
Namely, we use GM duality to produce a potential function that establishes weak reversal of the attracting posets between the mirror partners.
The relative GM identity gives the $\hbar=0$ part of the
Appell--Humbert quadratic-form comparison for elliptic stable envelopes.
Lastly, we show a combinatorial application: we interpret the core-quotient bijection for partitions as the 3d mirror correspondence on fixed points for certain bow varieties, and derive a per-color content formula in terms of the core and the quotient.
\end{abstract}

\maketitle

\section{Introduction}
\label{sec:intro}

\subsection{3d mirror symmetry and its shadows}
Physics predicts that certain holomorphic symplectic varieties with torus
actions come in pairs, called 3d mirror (or symplectic dual) pairs,
\[
 \X\ \circlearrowleft\ \T=\A\times\C^*_\h
 \qquad\longleftrightarrow\qquad
 \X^\dual\ \circlearrowleft\ \T^\dual=\A^\dual\times\C^*_{\h^\dual},
\]
where $\A$, $\A^\dual$ preserve and $\C^*_\h$, $\C^*_{\h^\dual}$ scale the
symplectic forms.  The mathematical content of the prediction consists of 
related statements at several levels, see e.g.\ \cite{WY23}.
Here is a(n incomplete) list: 
a bijection $f\leftrightarrow f^\dual$ between torus fixed points, a reversal of the
partial orders on the fixed points induced by attracting sets, 
an identification of the fixed-point restrictions of elliptic stable
envelopes \cite{AO21} after an exchange of equivariant and K\"ahler variables
(the \emph{stable-envelope mirror identity}) \cite{RSVZ22,BottaRimanyi}, the same identification for the
Appell--Humbert quadratic forms of the stable envelopes (the
\emph{quadratic-form mirror identity}), identifications of vertex functions \cite{BD25},
and equivalences of certain categories.

There is one more classical entry in this list, which lives on the level of degree two
equivariant cohomology.  For a torus $T$ acting on $\X$ with
isolated fixed points, the restriction $H^2_T(\X)\to H^2_T(p)$ to a fixed
point $p$ has a kernel $K_p$, and the collection of these kernels
is the degree-two part of the picture of $\operatorname{Spec}H^*_T(\X)$ studied by
Goresky and MacPherson \cite{GM10}.  They observed that the two arrangements are dual for $T^*(G/P)$ and
its dual slice, and for Gale dual pairs of hypertoric varieties.
More precisely, there is a perfect pairing between $H^2_T(\X)$ and $H^2_{T^\dual}(\X^\dual)$, under
which the equivariant parameters of one side are orthogonal to those of the
other, and $K_p^\perp=K_{p^\dual}$ for every fixed point.  Braden, Licata,
Proudfoot and Webster list this \emph{Goresky--MacPherson duality} among the
expected properties of a symplectic dual pair \cite[\S10]{BLPW16}.

\subsection{Bow varieties of affine type A} 
Bow varieties of affine type $A$ \cite{cherkis2011instantons,nakajima2017cherkis} are holomorphic
symplectic varieties attached to circular \emph{brane diagrams}.  Such a
diagram is a circle decorated with two types of five-branes, D5 and NS5,
with nonnegative integers on the arcs between them.  They generalize simultaneously the bow
varieties of type $A$ and the Nakajima quiver varieties of
affine type~$A$, including Hilbert schemes of points of $\C^2$, Hilbert schemes of
cyclic quotient resolutions of $\C^2$, and moduli spaces of framed torsion-free sheaves over $\C^2$.
A brane diagram $D$ has a combinatorial mirror $D^\dual$, obtained by
exchanging the two brane types.
\subsection{Results}
\label{subsec:results}
\emph{Standing convention.} Every circular brane diagram considered in this
paper has at least one five-brane of each type.  Denote such a diagram by $D$,
with D5 branes $a_1,\dots,a_n$ and NS5 branes $z_1,\dots,z_m$. Denote $P_D$ the finite
set of torus fixed points of $\X(D)$ and let $f\mapsto f^\dual$ be the mirror
bijection $P_D\to P_{D^\dual}$.  The symbols $a_i$, $b=t_1-t_2$ also denote
the equivariant parameters of the symplectic torus $\A_D$, and
$z_j$, $\kappa$ the K\"ahler parameters (Section~\ref{sec:recollections}).

Consider two finite-dimensional $\Q$-vector spaces carrying the data of these degree-two parameters.
\begin{itemize}
\item $E_D$, the character space of $\A_D$ modulo the framing diagonal, spanned
 by the $a_i-a_{i'}$ and $b$.  Its dimension is $n$.
\item $K_D$, the space of \emph{determinant classes}, with basis
 $\ell_1,\dots,\ell_m$.  Here $\ell_j$ stands for $\det\xi_{I_j}$, the
 determinant of the tautological bundle of the D3 interval $I_j$ just to the
 right of $z_j$.  Its dimension is $m$.
\end{itemize}
The mirror parameter map $\varphi_D$ (c.f. \eqref{eq:mirror-map}):
\[
 \begin{aligned}
 \varphi_D(a_j^\dual)=z_j-L_j\kappa,\quad \varphi_D(z_i^\dual)=a_i+R_i t_1,\\
 \varphi_D(t_1^\dual)=-\kappa,\quad\varphi_D(t_2^\dual)=\kappa-\h,\quad
 \varphi_D(\kappa^\dual)=-t_1,
 \end{aligned}
\] exchanging the equivariant and K\"ahler
parameters of $D$ and $D^\dual$, induces isomorphisms
\[
 \eta_D\colon K_{D^\dual}\xrightarrow{\ \cong\ }E_D^*,
 \qquad
 \eta_{D^\dual}\colon K_D\xrightarrow{\ \cong\ }E_{D^\dual}^*,
\]
and the \emph{GM pairings} (Definition~\ref{def:eta}).  

\smallskip

For $\lambda\in K_D$ and $f,g\in P_D$, let
$\delta^D_{f,g}(\lambda)\in E_D$ be the difference of the $\A_D$-weights of
the line bundle $L_D(\lambda)=\bigotimes_j(\det\xi_{I_j})^{\lambda_j}$ at $f$
and at $g$ for integral $\lambda$, extended $\Q$-linearly to all $K_D$.
\begin{mainthm}[Relative Goresky--MacPherson identity, Theorem~\ref{thm:A}]\label{mainthm:A}
For every circular brane diagram $D$ with $P_D\neq \varnothing$, all $f,g\in P_D$, all $\lambda\in K_D$
and $\lambda^\dual\in K_{D^\dual}$, we have
\[
 \bigl\langle\eta_D(\lambda^\dual),\ \delta^D_{f,g}(\lambda)\bigr\rangle
 +\bigl\langle\eta_{D^\dual}(\lambda),\ \delta^{D^\dual}_{f^\dual,g^\dual}(\lambda^\dual)\bigr\rangle
 =0 .
\]
\end{mainthm}

Denote $H_D=E_D\oplus K_D$ and let
$B\colon H_D\times H_{D^\dual}\to\Q$ be the pairing whose only nonzero blocks
are $\eta_D$ on $E_D\times K_{D^\dual}$ and $\eta_{D^\dual}$ on
$K_D\times E_{D^\dual}$, which looks like
\[
\begin{pmatrix}
0_{E_D \times E_D^!} & \ast_{K_D \times E_D^!} \\
\ast_{E_D \times K_D^!} & 0_{K_D \times K_D^!}
\end{pmatrix}.
\]

Assume $P_D\ne\varnothing$, fix a base point $f_0\in P_D$, and let
$V_f\subset H_D$ be the kernel of the restriction
$(\varepsilon,\lambda)\mapsto\varepsilon+\delta^D_{f,f_0}(\lambda)$, and
similarly $V^\dual_{f^\dual}\subset H_{D^\dual}$.

\begin{mainthm}[Goresky--MacPherson duality on the tautological model, Theorem~\ref{thm:C}]
\label{mainthm:C}
Assume $P_D\ne\varnothing$ and choose $f_0\in P_D$ to define the subspaces
above.  Then $B$ is perfect, $E_D\perp E_{D^\dual}$, $K_D\perp K_{D^\dual}$, and for every
$f\in P_D$ we have
\[
 (V_f)^{\perp_B}=V^\dual_{f^\dual}.
\]
\end{mainthm}

This is Goresky--MacPherson duality as in the literature, except that
$H^2_{\A}(\X)$ is replaced by its \emph{tautological model} $E_D\oplus K_D$,
in which the K\"ahler directions are the formal determinant classes rather
than $H^2(\X)$.  We regard this model as the natural habitat of GM duality for bow varieties.

To claim Goresky--MacPherson duality in the exact sense of \cite{GM10,BLPW16} we need some favorable conditions to hold.
Suppose that for both $\X$ and $\X^\dual$ the first Chern classes of the chosen tautological determinant line bundles form a basis of $H^2(\ ;\Q)$. Suppose also that the degree-two equivariant forgetful sequences\eqref{sequence} are exact for the symplectic tori
(modulo their framing diagonals). Then Theorem~\ref{mainthm:C} induces
Goresky--MacPherson duality on degree-two equivariant cohomology in the sense of \cite{GM10,BLPW16}, see Section~\ref{subsec:GM-comparison}.

\bigskip

The second theme of the paper concerns the Bruhat-like poset structure on bow varieties, and its (weak) reversal property for 3d dual pairs.

\begin{mainthm}[Weak order reversal, Theorem~\ref{thm:B}]
\label{mainthm:B}
Let $D$ be a circular brane diagram with $P_D\ne\varnothing$,
presented with a cut between two five-branes of opposite types.
Equip $\X(D)$ and $\X(D^\dual)$ with the chambers
in~\eqref{eq:chamber}, they assign attracting cells to fixed points, and hence an attracting poset on them. 
There is a function $H\colon P_D\to\mathbb R$
that is strictly increasing along the attracting order on $\X(D)$
and strictly decreasing along the attracting order on $\X(D^\dual)$
transported to $P_D$ by the mirror bijection. In particular, the pair
has weak order reversal in the sense of Definition~\ref{def:weak-or}.
\end{mainthm}

Conjecturally the exact (`strong') order reversal property also holds, but it is only known in type~A \cite{BFR24}. The weak one proved in this paper for type $\Ah$ can still serve as an ingredient in induction arguments on higher level 3dMS proofs \cite{BottaRimanyi-Ahat}.

\smallskip

Conjecturally, 3d mirror symmetry holds for the so-called elliptic stable envelope classes of Aganagic-Okounkov. This stable envelope 3dMS statement has a quadratic form `shadow', an identity for quadratic forms in the equivariant and K\"ahler variables. We prove that our GM duality is the $\h=0$ part of this quadratic form identity.

\smallskip

Since torus fixed points on bow varieties can be encoded with various combinatorial structures (tie diagrams, arrays of charged partitions) we deduce combinatorial consequences of the GM duality, along the lines of core-quotient bijections for partitions and per-color content identities.

\subsection{Organization}
Section~\ref{sec:recollections} recalls what we need about affine type $A$
bow varieties.
Section~\ref{sec:linear} sets up the spaces $E_D$, $K_D$ and the GM pairings.
Section~\ref{sec:GMidentity} states and proves the relative GM identity,
after working out the example of the Hilbert scheme of points of $\C^2$,
and records its absolute form.
Section~\ref{sec:GMduality} proves GM duality on the tautological model and
compares it with the formulation of \cite{GM10,BLPW16}, including the
base-point-free normalizations.
Section~\ref{sec:weakOR} uses the relative GM identity to construct a
height function and prove weak order reversal.
Section~\ref{sec:quadratic} identifies the relative GM identity with
the $\hbar=0$ part of the Appell--Humbert mirror comparison.
Section~\ref{sec:content} describes the mirror bijection using charged
partitions and derives the per-color content identity.

\subsection{Acknowledgments} 
The second author was supported by the NSF grants DMS-2200867. The authors are grateful to Tommaso Botta, \'Ad\'am Gyenge and Nok To Omega Tong for several helpful conversations. 
The authors provided the conceptual idea and framework of the theorems and proofs in this paper, but computer experiments, computations,  and independent audits of the arguments were produced with substantial help from AI systems (Anthropic's Claude and OpenAI's ChatGPT). The authors have checked all arguments and take full responsibility for the content.

\section{Recollections on bow varieties of affine type \texorpdfstring{$A$}{A}}
\label{sec:recollections}

We use the conventions of
\cite{rimanyi2020bow,GyR}.  
The reader familiar with these papers may skip to
Section~\ref{subsec:cells}, where the notation specific to this paper begins.

\subsection{Brane diagrams}
\label{subsec:diagrams}
A \emph{brane diagram of type $\Ah$} is an oriented circle carrying finitely
many marked points of two kinds, \emph{D5 branes} (drawn $\bs$) and
\emph{NS5 branes} (drawn $\fs$), together with a nonnegative integer $d_I$,
the \emph{multiplicity}, on each arc $I$ between consecutive five-branes.
The arcs are called \emph{D3 branes} or \emph{D3 intervals}.  A brane diagram of type $A$ is the
same on an oriented segment and it may be regarded as a circular diagram with
one arc of multiplicity zero.
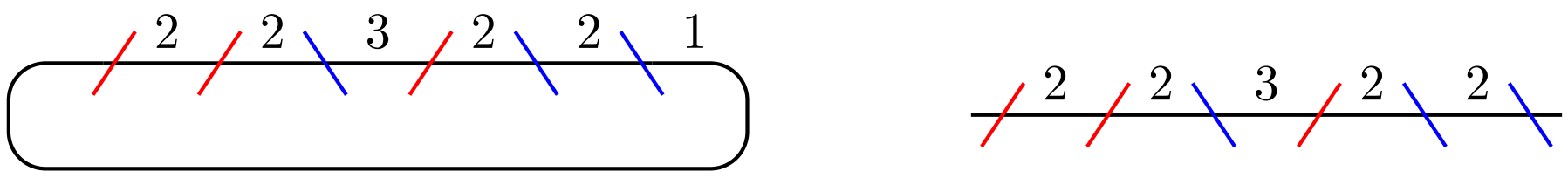
\begin{center}
\begin{tikzpicture}[scale=.8]
\draw[thick] (-0.1,0) -- (5.1,0);
\draw[red,thick] (-0.2,-0.3) -- (0.2,0.3);
\draw[red,thick] (0.8,-0.3) -- (1.2,0.3);
\draw[blue,thick] (1.8,0.3) -- (2.2,-0.3);
\draw[red,thick] (2.8,-0.3) -- (3.2,0.3);
\draw[blue,thick] (3.8,0.3) -- (4.2,-0.3);
\draw[blue,thick] (4.8,0.3) -- (5.2,-0.3);
\draw[rounded corners=8pt,thick]
  (5.1,0) -- (6,0) -- (6,-1) -- (-1,-1) -- (-1,0) -- (-0.1,0);
\draw (0.5,0.3) node {2};
\draw (1.5,0.3) node {2};
\draw (2.5,0.3) node {3};
\draw (3.5,0.3) node {2};
\draw (4.5,0.3) node {2};
\draw (5.5,0.3) node {1};
\end{tikzpicture}
\qquad\qquad
\begin{tikzpicture}[scale=.8,baseline=-12]
\draw[thick] (-0.3,0) -- (5.3,0);
\draw[red,thick] (-0.2,-0.3) -- (0.2,0.3);
\draw[red,thick] (0.8,-0.3) -- (1.2,0.3);
\draw[blue,thick] (1.8,0.3) -- (2.2,-0.3);
\draw[red,thick] (2.8,-0.3) -- (3.2,0.3);
\draw[blue,thick] (3.8,0.3) -- (4.2,-0.3);
\draw[blue,thick] (4.8,0.3) -- (5.2,-0.3);
\draw (0.5,0.3) node {2};
\draw (1.5,0.3) node {2};
\draw (2.5,0.3) node {3};
\draw (3.5,0.3) node {2};
\draw (4.5,0.3) node {2};
\end{tikzpicture}
\end{center}

To write formulas we \emph{cut} the circle at a point away from the
five-branes and unwrap it so that the orientation runs from left to right.
We draw the five-branes on top and the cut at the bottom.  The D5 branes are
then $a_1,\dots,a_n$ and the NS5 branes $z_1,\dots,z_m$, in left-to-right
order.  The same symbols denote the corresponding parameters.  So $a_i$ is an
equivariant parameter and $z_j$ is a K\"ahler parameter.  Relative to the cut
put
\begin{equation}\label{eq:R-L}
 R_i=\#\{j: z_j\text{ is to the right of }a_i\},
 \qquad
 L_j=\#\{i: a_i\text{ is to the left of }z_j\}.
\end{equation}
For a five-brane $v$ let $d_L(v)$, $d_R(v)$ be the multiplicities of the
adjacent intervals and define the \emph{local charges} (that are independent on the cut)
\begin{equation}\label{eq:lch}
 \lch(a_i)=d_L(a_i)-d_R(a_i),
 \qquad
 \lch(z_j)=d_R(z_j)-d_L(z_j).
\end{equation}

\emph{Standing convention on cuts.} Whenever chambers and attracting orders
are involved (Section~\ref{sec:weakOR}) we require the two five-branes
adjacent to the cut to be of opposite types (this requirement plays no role
in Sections~\ref{sec:GMidentity}--\ref{sec:GMduality}). The
requirement is symmetric under the mirror operation below.  

A \emph{Hanany--Witten transition} replaces a local configuration
$\ttt{$d_1$\fs $d_2$\bs $d_3$}$ by $\ttt{$d_1$\bs $d_1+d_3-d_2+1$\fs $d_3$}$ 
(in either direction).  A circular diagram is in \emph{normal form} with respect to a cut if all
NS5 branes are left of all D5 branes. Every diagram can be brought to normal form
by Hanany--Witten transitions. In \cite[Sect.~2.3]{GyR} the ambiguity of such normal forms is discussed.

\subsection{The mirror diagram}
The \emph{3d mirror} $D^\dual$ of $D$ is obtained by turning every D5 brane
into an NS5 brane and every NS5 brane into a D5 brane, keeping their positions,
the orientation, the multiplicities and the cut.  Thus $a_i$ becomes an NS5
brane $z^\dual_i$ of $D^\dual$ and $z_j$ becomes a D5 brane $a^\dual_j$.
The diagram $D^\dual$ has $m$ D5 and $n$ NS5 branes, and
\begin{equation}\label{eq:RL-dual}
 L^\dual_i=m-R_i,\qquad R^\dual_j=n-L_j,
 \qquad
 \lch(z^\dual_i)=-\lch(a_i),\qquad \lch(a^\dual_j)=-\lch(z_j).
\end{equation}
Pictorially, $D^\dual$ is the reflection of the picture of $D$ across the
circle.  The circle is fixed pointwise, inside and outside are exchanged, and
each little brane segment is reflected across the tangent line.  Clearly we have
$(D^\dual)^\dual=D$.

\subsection{Bow varieties (sketch)}
\label{subsec:bow}
Attached to the brane diagram $D$ there is a smooth holomorphic symplectic variety
$\X(D)$, the \emph{bow variety} \cite{cherkis2011instantons,nakajima2017cherkis}.

The bow variety is constructed from a space of linear maps between vector spaces
$V_I=\C^{d_I}$ and framing lines $\C_i$ (one per D5 brane).  Let
$\mathcal R_D\subset\mu_D^{-1}(0)$ be the locus satisfying the
conditions (S1), (S2) of \cite[\S2.2]{nakajima2017cherkis}.  This locus is affine,  and it is an open subset of the locus $\mathcal R_D\subset\mu_D^{-1}(0)$\cite[\S2.2(c), Proposition~3.2]{nakajima2017cherkis}.
Then we imposes the additional GIT stability condition for
$G_D=\prod_I GL(V_I)$ and the character
$\chi_0(g)=\prod_I\det(g_I)$, as in \cite[\S2.2]{GyR} and define
$\X(D)=\mathcal R_D/\!/_{\chi_0}G_D$. The torus
\begin{equation}\label{eq:torus}
 \T_D=\C^*_{t_1}\times\C^*_{t_2}\times\prod_{i=1}^n\C^*_{a_i}
\end{equation}
acts on $\X(D)$ as follows: $\C^*_{a_i}$ rescales the framing
line $\C_i$, and $t_1,t_2$ rescale the linear maps.  Put
\[
 \h=t_1t_2,\qquad b=t_1/t_2 .
\]
The character $\h$ scales the symplectic form and its kernel
$\A_D\subset\T_D$ is the \emph{symplectic torus}.  The diagonal subtorus of
$\prod_i\C^*_{a_i}$ acts trivially. The fixed point set $\X(D)^{\T_D}$ is
finite and equals to $\X(D)^{\A_D}$.

Each D3 interval $I$ carries a $\T_D$-equivariant \emph{tautological bundle}
$\xi_I$ of rank $d_I$, induced from the standard representation of $GL(V_I)$.
Its Chern roots are $x_{I,1},\dots,x_{I,d_I}$, and we will use the notation 
$X_I=\{x_{I,1},\dots,x_{I,d_I}\}$.

Under a choice of cut, for the $j-$th NS5 brane we assign a `K\"ahler parameter' $z_j$. Moreover we introduce a 
\emph{normalized K\"ahler period} $\kappa$ for the diagram $D$, with the intuitive meaning that the `K\"ahler monodromy around the circle' is $\zeta=\kappa^n$.  These are formal parameters, not characters of~$\T_D$.

\emph{Additive and multiplicative notation.}  All of $a_i,z_j,t_1,t_2,\h,b,
\kappa,\zeta$ appear both multiplicatively (as characters, in $K$-theory and
elliptic cohomology) and additively (as logarithms of multiplicative characters, in cohomology and in
the linear algebra of this paper).  In additive notation $\h=t_1+t_2$,
$b=t_1-t_2$, $\zeta=n\kappa$.  From Section~\ref{sec:linear} on, all
parameters are additive, we will make exceptional mentioning if we shift to the multiplicative convention.

If two diagrams differ by a Hanany--Witten transition, the corresponding bow
varieties are isomorphic.  For the direction moving the D5 to the right
past the NS5, the isomorphism is equivariant after the reparametrization
$\widetilde a=a+t_1$ of the parameter of the moving D5 brane (use the inverse substitution for
the reverse move). Denoting the three tautological bundles before the transition by $\xi_1,\xi_2,\xi_3$ and the new middle bundle by $\widetilde{\xi}_2$, the 
tautological bundles are related by the exact sequence
\begin{equation}\label{eq:HW-taut}
 0\to t_2^{-1}\xi_2\to t_2^{-1}\xi_1\oplus\xi_3\oplus\C_{\widetilde a}
 \to\widetilde\xi_2\to0 .
\end{equation}

\subsection{Fixed points and tie diagrams}
\label{subsec:ties}
\begin{figure}[ht]
\[
\begin{tikzpicture}[scale=0.7]
    \draw[rounded corners=10mm, thick] (0,5) -- (10,5) -- (10,10) --(0,10) -- cycle;
    \draw[thick,red] (2,9.6)--(2.4,10.4);
    \draw[thick,red] (3,9.6)--(3.4,10.4);
    \draw[thick,blue] (6.4,9.6)--(6,10.4);
    \draw[thick,blue] (7.4,9.6)--(7,10.4);
    \draw[thick,blue] (8.4,9.6)--(8,10.4);
\draw[thick,dotted] (6,10.4) to [out=150,in=30] (3.4,10.4);
\draw[thick,dotted] (7,10.4) to [out=155,in=30] (2.4,10.4);
\draw[thick,dotted,rounded corners=12mm]
  (7,10.4) to [out=150,in=30] (-.5,10) -- (-.5,4.5)--
  (10.5,4.5) -- (10.5,10.5) to [out=170,in=35] (3.4,10.4);
\draw[thick,dotted,rounded corners=9mm]
  (8.4,9.6) -- (9.5,8.5) -- (9.5,5.5) -- (0.5,5.5) --
  (0.5,9) -- (3,9.6);
\draw[thick,dotted,rounded corners=9mm]
  (8.4,9.6) -- (9,8.5) -- (9,6) -- (1,6) -- (1,8.8) --
  (8.5,8.8) -- (8.5,6.5) -- (1.5,6.5) -- (1.5,8.8) -- (3,9.6);
\draw (2,10.4) node {$z_1$};
\draw (3.3,9.6) node {$z_2$};
\draw (6,9.6) node {$a_1$};
\draw (7,9.6) node {$a_2$};
\draw (8,9.6) node {$a_3$};
\end{tikzpicture}
\]
\[
\begin{tikzpicture}[scale=.55]
\draw[thick]  (0.5,0) -- (16,0);
\draw[thick,red] (1,-.4)--(1.5,.4);
\draw[thick,red] (1.5,-.4)--(2,.4);
\draw[thick,blue] (3.5,-.4)--(3,.4);
\draw[thick,blue] (4,-.4)--(3.5,.4);
\draw[thick,blue] (4.5,-.4)--(4,.4);
\draw[thick,red] (5.5,-.4)--(6,.4);
\draw[thick,red] (6,-.4)--(6.5,.4);
\draw[thick,blue] (8,-.4)--(7.5,.4);
\draw[thick,blue] (8.5,-.4)--(8,.4);
\draw[thick,blue] (9,-.4)--(8.5,.4);
\draw[thick,red] (10,-.4)--(10.5,.4);
\draw[thick,red] (10.5,-.4)--(11,.4);
\draw[thick,blue] (12.5,-.4)--(12,.4);
\draw[thick,blue] (13,-.4)--(12.5,.4);
\draw[thick,blue] (13.5,-.4)--(13,.4);
\draw[thick,red] (14.5,-.4)--(15,.4);
\draw[thick,red] (15,-.4)--(15.5,.4);
\draw[thick,dotted] (2,0.4) to [out=30,in=150] (8,0.4);
\draw[thick,dotted] (6,0.4) to [out=30,in=150] (8,0.4);
\draw[thick,dotted] (6.5,0.4) to [out=30,in=150] (7.5,0.4);
\draw[thick,dotted] (9,-0.4) to [out=-30,in=-150] (10.5,-0.4);
\draw[thick,dotted] (9,-0.4) to [out=-30,in=-150] (15,-0.4);
\draw (8,-.4)--(7.1,-1.1) node[below] {$a_1$};
\draw (8.5,-.4)--(8,-1.1) node[below] {$a_2$};
\draw (9,-.4)--(8.9,-1.1) node[below] {$a_3$};
\end{tikzpicture}
\]
\caption{A tie diagram of a circular brane diagram
($\ttt{$4$\fs $5$\fs $5$\bs $4$\bs $2$\bs}$ on the circle)
and the picture over the universal cover. The lifted endpoints distinguish
ties with different winding numbers.}
\label{fig:tiediagrams}
\end{figure}
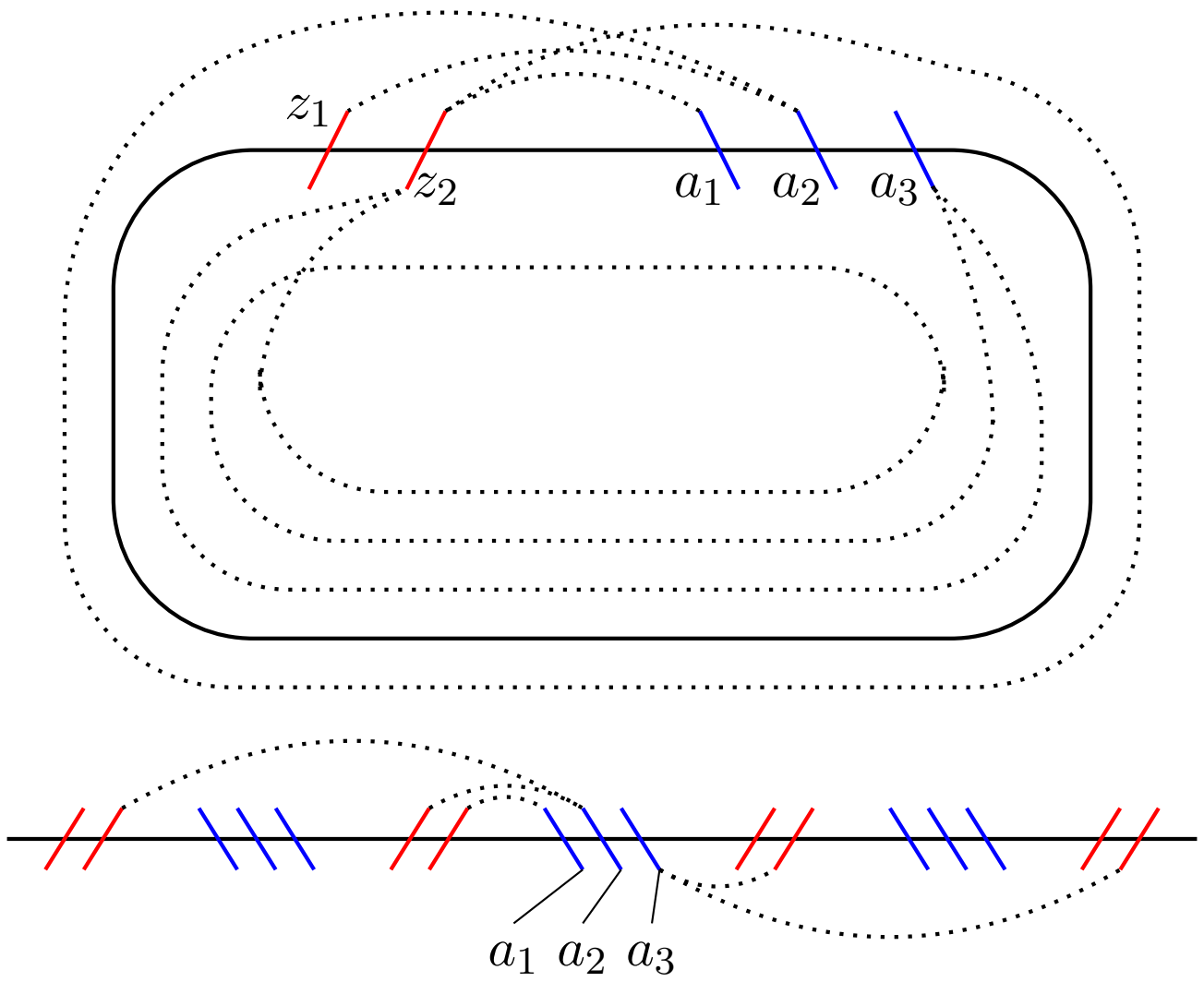

A \emph{tie} connects a D5 brane and an NS5 brane. In a circular diagram it may wind around the circle, that is, it also remembers how many times it winds around in which direction.  On the universal cover of the circle, a tie is a segment between a lift of its D5 endpoint and a lift of its NS5 endpoint, and it \emph{covers} the lifted D3 intervals in between.  A
\emph{tie diagram} is a set of ties, such that every lifted D3 interval $I$ is covered by exactly $d_I$
ties. Let us emphasize that the diagram is a {\em set}, not a multiset: each triple (D5, NS5, winding) corresponds to at most one tie. We have a bijection (see \cite{nakajima2017cherkis,rimanyi2020bow,GyR}) 
\begin{equation}\label{eq:fixed-points-ties} 
 P_D:=\X(D)^{\T_D}\ \longleftrightarrow\ \{\text{tie diagrams of }D\}.
\end{equation}
Figure~\ref{fig:tiediagrams} shows a tie diagram of a circular diagram and a
window of its universal cover.

If $D$ and $\widetilde D$ differ by a Hanany--Witten transition, the
isomorphism $\X(D)\cong\X(\widetilde D)$ induces a bijection of fixed points, the bijection is explicit on tie diagrams: it is a local move carrying the ties attached to the two
transposed branes across each other, see \cite[Fig.~5.]{rimanyi2020bow}.

\subsection{Mirror fixed points}
\label{subsec:mirror-fixed-points}
\begin{definition} 
\label{def:mirror-fixed-point}
Let $f\in P_D$ be given by a tie diagram. Reflecting the tie diagram across the circle (i.e. tie drawn inside the circle
moves outside and vice versa) exchanges the D5 and NS5 branes graphically. The fixed point $f^\dual$ represented by this reflected tie diagram will be called the \emph{mirror fixed point} in $P_{D^!}$ of mirror bow variety $\X^!(D)$.
\end{definition}

Since the reflection preserves all covering numbers, $f^\dual \in P_{D^\dual}$. Clearly $(f^\dual)^\dual=f$.

\subsection{Cells and the signed tie indicator}
\label{subsec:cells}
Let $D$ be in normal form with the cut between $a_n$ and $z_1$, so that one
period of the universal cover reads:
\[
 z_1^{(r)}\cdots z_m^{(r)}\ a_1^{(r)}\cdots a_n^{(r)},\qquad r\in\Z.
\]
Here the superscript labels the various lifts and the deck transformation is
$r\mapsto r+1$.  Thus $z_j^{(s)}$ lies to the left of $a_i^{(r)}$ iff
$s\le r$.  Every tie has a unique lift whose D5 endpoint is $a_i^{(0)}$ for
some $i$.  Its NS5 endpoint is then some $z_j^{(p)}$, and we call the tie the
\emph{cell} $(p,i,j)\in\Z\times\{1,\dots,n\}\times\{1,\dots,m\}$. (The word `cell' refers to the entry position in the 3d matrix $\Z\times\{1,\dots,n\}\times\{1,\dots,m\}$.)

The tie is
a \emph{left tie} if $p\le0$ (the NS5 endpoint is to the left of the D5
endpoint) and a \emph{right tie} if $p\ge1$, cf.\ \cite[App.~A]{nakajima2018towards}, \cite[App~A]{rimanyi2020bow}.

For $f\in P_D$ define the \emph{signed tie indicator}
\begin{equation}\label{eq:N}
 \nu_f(p,i,j)=
 \begin{cases}
 +1&\text{if the cell $(p,i,j)$ is a (left) tie of $f$ and $p\le0$},\\
 -1&\text{if the cell $(p,i,j)$ is a (right) tie of $f$ and $p\ge1$},\\
 0&\text{otherwise.}
 \end{cases}
\end{equation}
It clearly has finite support and determines $f$.  In the language of
\cite{GyR}, $\nu_f=M_f-M^{\mathrm{vac}}$ is the difference between the affine
Maya diagram of $f$ and the vacuum diagram, with the half-integer block
coordinate $q$ replaced by $p=q+\frac12$.

\begin{proposition}[\cite{GyR}]\label{prop:invariants}
For every $f\in P_D$ we have
\begin{equation}\label{eq:invariants}
\begin{gathered}
 -\sum_{p,j}\nu_f(p,i,j)=-\lch(a_i)\quad(1\le i\le n),\qquad
 -\sum_{p,i}\nu_f(p,i,j)=-\lch(z_j)\quad(1\le j\le m),\\
 \sum_{p,i,j}(1-p)\,\nu_f(p,i,j)=d_{I_m},
\end{gathered}
\end{equation}
where $I_m$ is the D3 interval between $z_m$ and $a_1$.
\end{proposition}
In the language of Maya diagrams of \cite{GyR}, the three quantities are the row charges $e_i$, the column charges $f_j$ and
the level.  For $f,g\in P_D$ we have
\begin{equation}\label{eq:conservation}
 \sum_{p,j}(\nu_g-\nu_f)(p,i,j)=0\ \ \forall i,\qquad
 \sum_{p,i}(\nu_g-\nu_f)(p,i,j)=0\ \ \forall j,\qquad
 \sum_{p,i,j}p\,(\nu_g-\nu_f)(p,i,j)=0 .
\end{equation}

\begin{lemma}\label{lem:mirror}
Let $D$ be in normal form as above.  Present $D^\dual$ in its own normal form
$z^\dual_1\cdots z^\dual_n\,a^\dual_1\cdots a^\dual_m$ by choosing the lifts
$z^{\dual(0)}_i$ at the position of $a_i^{(0)}$ and $a^{\dual(0)}_j$ at the
position of $z_j^{(1)}$. Namely, if the 5 branes in D look
\[
 \ldots\underbrace{z_1^{(0)}\cdots z_m^{(0)}\ a_1^{(0)}\cdots a_n^{(0)}}z_1^{(1)}\cdots z_m^{(1)}\ldots,
\]
then its mirror looks 
\[
\ldots a_1^{!(0)}\cdots a_m^{!(0)}\ \underbrace{z_1^{!(0)}\cdots z_n^{!(0)}a_1^{!(1)}\cdots a_m^{!(1)}}\ldots .
\]
The underbrace indicates the normal form interval.

Then for all $f\in P_D$
\begin{equation}\label{eq:N-mirror}
 \nu_{f^\dual}(1-p,\,j,\,i)=-\nu_f(p,i,j).
\end{equation}
\end{lemma}
\begin{proof}
The tie $(p,i,j)$ of $D$ connects $a_i^{(0)}$ and $z_j^{(p)}$.  Its
reflection connects $z^{\dual(0)}_i$ and the D5 brane of $D^\dual$ at the
position of $z_j^{(p)}$, which is $a^{\dual(p-1)}_j$.  Translating by $1-p$
periods, this is the tie between $a^{\dual(0)}_j$ and $z^{\dual(1-p)}_i$,
i.e.\ the cell $(1-p,j,i)$ of $D^\dual$.  A left tie of $D$ ($p\le0$) has its
NS5 endpoint to the left of its D5 endpoint.  After reflection the roles of the
endpoints are exchanged, so it is a right tie of $D^\dual$ ($1-p\ge1$), and
vice versa.  This explain the sign.
\end{proof}
\subsection{Weights of tautological bundles at fixed points}
\label{subsec:decoration}
At a fixed point $f$ the bundle $\xi_I$ splits into $\T_D$-weight lines, one
for each tie portion covering $I$. These weights are read off from the tie
diagram by the following rule (the affine version of the ``decorated butterfly'' rule of \cite[\S4.6]{rimanyi2020bow}). 

Fix a D5 brane $a_i$ and a lift.  Order the ties ending at $a_i$
from the left from longest to shortest. The portion just left of $a_i$ on the
tie of rank $q\ge0$ has weight $a_i-q\h$, and moving along the tie to the
left across an NS5 brane subtracts $t_1$. 

Here is an example with $3$ NS5 branes and some D5 branes $\{a_k\}$, the weight under ties ending at $j-th$ D5 branes $a_j$ from the left can be read off as the following picture.
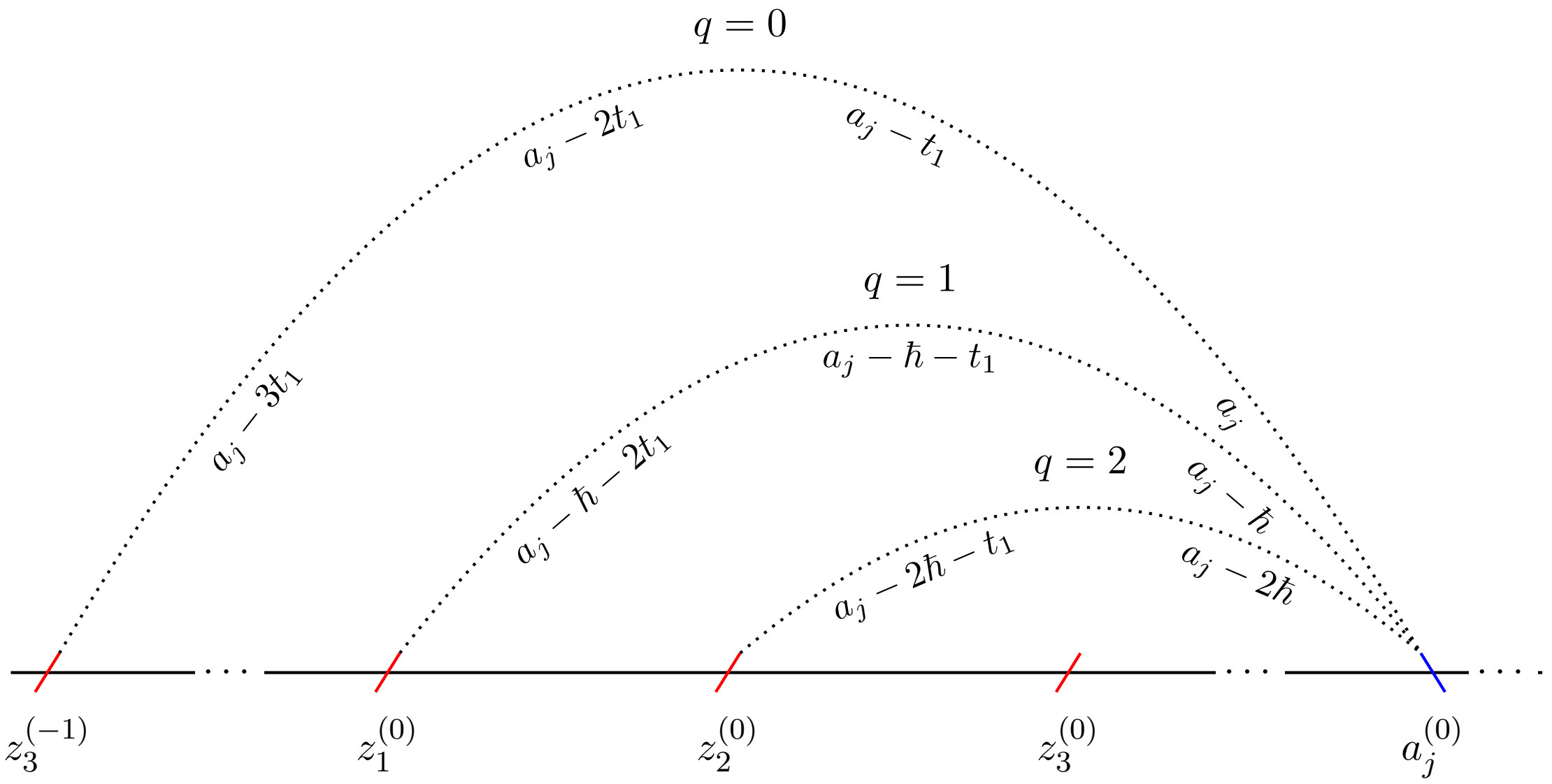
\[
\begin{tikzpicture}[
  scale=.46,
  weight/.style={
    sloped,
    below=3pt,
    font=\small,
    inner sep=1pt
  }
]

  \draw[thick,dotted]
    (1,.4)
    .. controls (10.333333,16.4) and (19.666667,16.4) ..
    node[pos=.125,weight] {$a_j-3t_1$}
    node[pos=.375,weight] {$a_j-2t_1$}
    node[pos=.625,weight] {$a_j-t_1$}
    node[pos=.875,weight] {$a_j$}
    node[pos=.5,above=3pt] {$q=0$}
    (29,.4);

  \draw[thick,dotted]
    (8,.4)
    .. controls (15,9.4) and (22,9.4) ..
    node[pos=.166667,weight] {$a_j-\hbar-2t_1$}
    node[pos=.5,weight] {$a_j-\hbar-t_1$}
    node[pos=.833333,weight] {$a_j-\hbar$}
    node[pos=.5,above=3pt] {$q=1$}
    (29,.4);

  \draw[thick,dotted]
    (15,.4)
    .. controls (19.666667,4.4) and (24.333333,4.4) ..
    node[pos=.25,weight] {$a_j-2\hbar-t_1$}
    node[pos=.75,weight] {$a_j-2\hbar$}
    node[pos=.5,above=3pt] {$q=2$}
    (29,.4);

  \draw[thick] (0,0) -- (31.5,0);

  \foreach \x in {1,8,15,22}{
    \draw[thick,red] (\x-.5,-.4) -- (\x,.4);
  }

  \draw[thick,blue] (29.5,-.4) -- (29,.4);

  \node[fill=white,inner sep=2pt] at (4.5,0)  {$\cdots$};
  \node[fill=white,inner sep=2pt] at (25.5,0) {$\cdots$};
  \node[fill=white,inner sep=2pt] at (30.7,0) {$\cdots$};

  \node[below] at (.75,-.6)   {$z_3^{(-1)}$};
  \node[below] at (7.75,-.6)  {$z_1^{(0)}$};
  \node[below] at (14.75,-.6) {$z_2^{(0)}$};
  \node[below] at (21.75,-.6) {$z_3^{(0)}$};
  \node[below] at (29.25,-.6) {$a_j^{(0)}$};
\end{tikzpicture}
\]
If $\ell$ ties end at $a_i$ from
the left, order the ties ending at $a_i$ from the right from shortest to
longest.  The portion just right of $a_i$ on the tie of rank $q\ge0$ has weight
$a_i+(1-\ell+q)\h$, and moving to the right across an NS5 brane subtracts
$t_2$.  Ties passing through $a_i$ keep their weight.  The multiset of weights
of the portions covering a lifted interval $I$ is $X_I|_f$.

We record two consequences.

\begin{lemma}\label{lem:jump}
For every D5 brane $a_i$ and every fixed point,
\[
 \sum_{x\in X_{I_R(a_i)}}x-\sum_{x\in X_{I_L(a_i)}}x
 =-\lch(a_i)\,a_i+\binom{-\lch(a_i)+1}{2}\h ,
\]
independently of the fixed point.  Consequently, if $I$ and $I'$ are D3
intervals in between the same two consecutive NS5 branes, and denote the D5 branes in between by $a_{m+1}...,a_{m+r}$, the weights of
$\det\xi_I$ and $\det\xi_{I'}$ at every fixed point differ by \[\sum_{i=1}^r \bigr(\sum_{x\in X_{I_R(a_{m+i})}}x-\sum_{x\in X_{I_L(a_{m+i})}}x\bigr)=\sum_i^r-\lch(a_{m+i})\,a_{m+i}+\binom{-\lch(a_{m+i})+1}{2}\h\] a linear
polynomial in $a,\h$, independent of fixed point.
\end{lemma}

Modulo $\h$ a cleaner identity is
\begin{equation}\label{eq:D5-virtual-alphabet}
 \overline X_{I_L(a_i)}-\overline X_{I_R(a_i)}=\lch(a_i)\{a_i\}.
\end{equation}
This follows directly from the decoration rule.  Passing ties cancel, while
all portions ending at the D5 have weight $a_i$ modulo $\h$.

\begin{lemma}[Tie-intrinsic weights modulo $\h$]\label{lem:intrinsic}
Let $t$ be a tie with D5 endpoint $a_i$ and let $I$ be a lifted D3 interval
covered by $t$, separated from the D5 endpoint by $k$ lifted NS5 branes.
Modulo $\h$, the weight of the portion of $t$ over $I$ is $a_i-\frac k2b$ if
$I$ is to the left of the D5 endpoint and $a_i+\frac k2b$ if it is to the
right.  In particular it depends only on $t$ and $I$, not on the fixed point.
\end{lemma}
\begin{proof}
The weight is $a_i-q\h-kt_1$ on the left and $a_i+(1-\ell+q)\h-kt_2$ on
the right, where $q$ and $\ell$ depend on the fixed point. Also, modulo $\h$ we have
$t_1\equiv\frac b2$ and $t_2\equiv-\frac b2$.
\end{proof}

\subsection{The mirror parameter map}
\label{subsec:mirror-map}
For a cut circular diagram $D$, we define the \emph{mirror parameter map} $\varphi_D$
to be the additive substitution
\begin{equation}\label{eq:mirror-map}
{\renewcommand{\arraystretch}{1.3}
\begin{array}{ll}
 \varphi_D(a_j^\dual)=z_j-L_j\kappa,  &\qquad  \varphi_D(z_i^\dual)=a_i+R_it_1,  \\
 \varphi_D(t_1^\dual)=-\kappa,        &\qquad \varphi_D(\kappa^\dual)=-t_1,\\
 \varphi_D(t_2^\dual)=\kappa-\h.     &
 \end{array}}
\end{equation}
Calculation gives $\varphi_D(\h^\dual)=-\h$, $\varphi_D(b^\dual)=\h-2\kappa$ and
$\varphi_D(\zeta^\dual)=\varphi_D(m\kappa^\dual)=-mt_1$, that is, $\varphi_D$ is an involution up to
deck transformations. Namely, by \eqref{eq:RL-dual} we have
\begin{equation}\label{eq:mirror-map-square}
 (\varphi_D\circ\varphi_{D^\dual})(a_i)=a_i+mt_1,\qquad
 (\varphi_D\circ\varphi_{D^\dual})(z_j)=z_j-\zeta,
\end{equation}
and $t_1,t_2,\kappa$ are fixed.  The combinations
\begin{equation}\label{eq:A}
 A_i=a_i+R_it_1=\varphi_D(z^\dual_i),
 \qquad
 A^\dual_j=a^\dual_j+R^\dual_jt^\dual_1
\end{equation}
are invariant under Hanany--Witten transitions.
\subsection{Change of cut}
Let \(r=\sum_i\chi_i\) and \(s=\sum_j\eta_j\) be the numbers of D5 and NS5 lifts crossing the cut, respectively. We have: \begin{align}
 a_i^c&=a_i^0+m\chi_i t_1,&
 z_j^c&=z_j^0+\eta_j\zeta,
 \label{eq:cut-coordinates}\\
 L_j^c&=L_j^0-r+n\eta_j,&
 R_i^c&=R_i^0+s-m\chi_i.
 \label{eq:cut-counts}
\end{align}
Here the quantities with superscript $c$ represent quantities after change of cut, while quantities with subscript $0$ represents quantities before the change. Then $A_i$ and $A^!_j$ change only by common
deck shifts under changes of cut 
\[
A_i^{c}
=(a_i^{0}+m\chi_i t_1)+(R_i^{0}+s-m\chi_i)t_1
=A_i^{0}+s\,t_1
=\varphi_{D^{c}}(z_i^{!,c}).
\]
On the dual diagram $D^!$ we have
\[
A_j^{!,c}
=(a_j^{!,0}+n\eta_j t_1^!)
 +(R_j^{!,0}+r-n\eta_j)t_1^!
=A_j^{!,0}+r\,t_1^!.
\]

\section{The linear algebra of parameters}
\label{sec:linear}

From now on all parameters are additive, and $D$ is presented with a cut.
We retain the standing convention of Section~\ref{subsec:results}.

\subsection{Equivariant characters}
\label{subsec:E}
On the symplectic torus $\A_D=\ker\h$ one has $t_1\equiv-t_2\equiv\frac b2$.
For an additive character $w=\sum_ic_ia_i+pt_1+qt_2$ of $\T_D$, write
\begin{equation}\label{eq:reduction}
 \ol w=\sum_ic_ia_i+\frac{p-q}2\,b\ \in X^*(\A_D)_\Q
\end{equation}
for the reduction of $w$ modulo $\h$.  Since the diagonal of $\prod_i\C^*_{a_i}$
acts trivially, every tangent weight at a fixed point, as well as every difference of
weights of a tautological line bundle at two fixed points, lie in the
\emph{reduced character space}
\[
 E_D=\operatorname{span}_\Q\{a_i-a_{i'},\ b\}\subset X^*(\A_D)_\Q,
 \qquad \dim E_D=n .
\]
Elements of $E_D^*$ are rational cocharacters of $\A_D$ modulo the diagonal.
For $\sigma\in E_D^*$, choose a lift and denote $\sigma_{a_i}=\pair\sigma{a_i}$, which is well defined up
to a common constant for all $i$. Define also $\sigma_{t_1}=\pair\sigma{t_1}=\frac12\pair\sigma b$. We have 
$\sigma_{t_2}=-\sigma_{t_1}$ and
\begin{equation}\label{eq:tangent-pairing}
 \pair{\sigma}{a_{i'}-a_i+pt_1+qt_2}
 =\sigma_{a_{i'}}-\sigma_{a_i}+(p-q)\sigma_{t_1}, \quad p-q\in\Z.
\end{equation}
Every tangent weight of $\X(D)$ at a fixed point has this form, and very importantly $p-q\ne0$ if $i=i'$. They cannot be zero because the fixed points are isolated
already for $\A_D$.

\subsection{K\"ahler characters and determinant classes}
\label{subsec:K}
Let $\widetilde z_j$ denote the K\"ahler coordinate of the chosen lift of
$z_j$, then the lift one period to the right has coordinate
$\widetilde z_j+\zeta$.  Put
\begin{equation}\label{eq:beta}
 \beta_j=\widetilde z_j-\widetilde z_{j+1}\ \ (1\le j<m),
 \qquad
 \beta_m=\widetilde z_m-\widetilde z_1-\zeta,
 \qquad
 \sum_{j=1}^m\beta_j=-\zeta=-n\kappa,
\end{equation}
the $m$ consecutive differences of lifted K\"ahler coordinates, and
\[
 Z_D=\operatorname{span}_\Q\{\beta_1,\dots,\beta_m\}
 =\operatorname{span}_\Q\{z_j-z_{j+1}\ (j<m),\ \kappa\},
 \qquad\dim Z_D=m .
\]
Let $I_j$ be the first D3 interval to the right of $z_j$ and let $\ell_j$ be
a formal symbol for the class of $\det\xi_{I_j}$.  The
\emph{determinant space} is
\[
 K_D=\bigoplus_{j=1}^m\Q\,\ell_j ,
\]
and we identify $Z_D=K_D^*$ by $\pair{\beta_j}{\ell_k}=\delta_{jk}$.  

We call a class
$\lambda=\sum_j\lambda_j\ell_j$ is \emph{positive} if all $\lambda_j>0$. For integral $\lambda$ we write
\[
 L_D(\lambda)=\bigotimes_{j=1}^m(\det\xi_{I_j})^{\lambda_j}
\]
with its tautological $\T_D$-equivariant structure.  By Lemma~\ref{lem:jump},
replacing $I_j$ by another interval in the same NS5 gap changes the weights of
$L_D(\lambda)$ at all fixed points by the same character.

For $\lambda\in K_D$ and $f\in P_D$ let
$\wt(\lambda|_f)=\sum_j\lambda_j\sum_{x\in X_{I_j}}x|_f\in X^*(\T_D)_\Q$ be
the weight of $L_D(\lambda)$ at~$f$.  

\begin{definition}
For $f,g\in P_D$ define the
\emph{relative restriction map}
\begin{equation}\label{eq:delta}
 \delta^D_{f,g} :  K_D\to E_D,
 \qquad
 \delta^D_{f,g}(\lambda)=\ol{\wt(\lambda|_f)-\wt(\lambda|_g)} .
\end{equation}
\end{definition}
The map lands in $E_D$ because of the following argument. Every tautological weight has $a$-degree one, and therefore $l_i|_f=n_i^fb+\sum_{j\geq 1} m_j^fa_j$ for some $m_j\geq 0$ such that $\sum m_j=d_{I_i}$. Then we have
\[
l_i|_f-l_i|_g=(n_0^f-n_0^g)b+\sum_j (m_j^f-m_j^g)a_j
\]
with $\sum_j (m_j^f-m_j^g)=0$. This is exactly the condition of being an element in $E_D$. 

The map $ \delta^D_{f,g}$  
satisfies
\[
\delta_{f,g}+\delta_{g,h}=\delta_{f,h}, \qquad \delta_{g,f}=-\delta_{f,g}.
\]

\subsection{The blocks of the mirror parameter map}
\begin{lemma}\label{lem:blocks}
Modulo $\h$, $\varphi_D$ maps $E_{D^\dual}$ into $Z_D$ and $Z_{D^\dual}$ into
$E_D$.  The induced maps
\[
 u_D\colon Z_{D^\dual}\to E_D,\qquad v_D\colon E_{D^\dual}\to Z_D
\]
are isomorphisms, given on bases by
\begin{align}
 u_D(\beta^\dual_i)&=\ol{A_i-A_{i+1}}\ \ (i<n),&
 u_D(\beta^\dual_n)&=\ol{A_n-A_1+mt_1},\label{eq:u}\\
 v_D(a^\dual_j-a^\dual_{j+1})&=\beta_j-(L_j-L_{j+1})\kappa\ \ (j<m),&
 v_D(b^\dual)&=-2\kappa .\label{eq:v}
\end{align}
\end{lemma}
\begin{proof}
Formulas \eqref{eq:u}--\eqref{eq:v} are \eqref{eq:mirror-map} read modulo
$\h$.  The wrap term uses $\varphi_D(\zeta^\dual)=-mt_1$.  The characters
$\ol{A_i-A_{i+1}}$ $(i<n)$ together with $\ol{t_1}=\frac b2$ span $E_D$, and
the $\beta_j$ $(j<m)$ together with $\kappa$ span $Z_D$.  In these bases both
maps are triangular with nonzero diagonal.
\end{proof}

\begin{definition}[The GM pairings]\label{def:eta}
Let $u_D^*\colon E_D^*\to Z_{D^\dual}^*=K_{D^\dual}$ and
$v_D^*\colon Z_D^*=K_D\to E_{D^\dual}^*$ be the transposes and we define
\[
 \eta_D=(u_D^*)^{-1}\colon K_{D^\dual}\xrightarrow{\ \cong\ }E_D^*,
 \qquad
 \eta_{D^\dual}=v_D^*\colon K_D\xrightarrow{\ \cong\ }E_{D^\dual}^*.
\]
\end{definition}
We observed in \eqref{eq:mirror-map-square} that $\varphi_{D^\dual}\circ\varphi_D$ is a deck
transformation, which acts on $E$ and $Z$ by common shifts, i.e.\ trivially.
Hence $v_D=u_{D^\dual}^{-1}$ and $\eta_{D^\dual}=(u_{D^\dual}^*)^{-1}$, and the two
definitions are symmetric under $D\leftrightarrow D^\dual$.

\begin{lemma}[The GM pairings in coordinates]\label{lem:eta}
For $\lambda^\dual=\sum_i\lambda^\dual_i\ell^\dual_i\in K_{D^\dual}$, denote
$T^\dual_\lambda=\sum_i\lambda^\dual_i$.  Then $\sigma=\eta_D(\lambda^\dual)$ is the unique element
of $E_D^*$ such that
\begin{equation}\label{eq:eta-char}
 \pair{\sigma}{A_i-A_{i+1}}=\lambda^\dual_i\ \ (1\le i<n),
 \qquad
 \pair{\sigma}{A_n-A_1+mt_1}=\lambda^\dual_n ,
\end{equation}
and consequently,
\begin{equation}\label{eq:eta-t1}
 \sigma_{t_1}=\langle \sigma, t_1\rangle=\frac{T^\dual_\lambda}m,
 \qquad
 \sigma_{a_{i+1}}-\sigma_{a_i}=\pair{\sigma}{a_{i+1}-a_i}=-\lambda^\dual_i+(R_i-R_{i+1})\frac{T^\dual_\lambda}m
 \quad(1\le i<n).
\end{equation}
Symmetrically, for $\lambda=\sum_j\lambda_j\ell_j\in K_D$ with
$T_\lambda=\sum_j\lambda_j$, the element $\sigma^\dual=\eta_{D^\dual}(\lambda)$
is characterized by $\pair{\sigma^\dual}{A^\dual_j-A^\dual_{j+1}}=\lambda_j$
$(j<m)$ and $\pair{\sigma^\dual}{A^\dual_m-A^\dual_1+nt^\dual_1}=\lambda_m$,
and $\sigma^\dual_{t^\dual_1}=T_\lambda/n$.
\end{lemma}
\begin{proof}
By definition $u_D^*(\sigma)=\sum_i\pair{\sigma}{u_D(\beta^\dual_i)}\ell^\dual_i$,
by adjunction, $\eta_D(\lambda^\dual)=\sigma$ iff $\pair{\sigma}{u_D(\beta^\dual_i)}=\lambda^\dual_i$ for
all $i$, which is \eqref{eq:eta-char} by \eqref{eq:u}.  

Summing the $n$ equations \eqref{eq:eta-char} gives $m\sigma_{t_1}=T^\dual_\lambda$.  

Since $a_i=A_i-R_it_1$, we have,
$\sigma_{a_{i+1}}-\sigma_{a_i}=-\pair\sigma{A_i-A_{i+1}}-(R_{i+1}-R_i)\sigma_{t_1}$,
This is exactly \eqref{eq:eta-t1}.  

The proof for the mirror side is essentially the same, hence we omit it.
\end{proof}

\section{The relative GM identity}
\label{sec:GMidentity}

\subsection{Statement}
\begin{theorem}[Relative GM identity]\label{thm:A}
Let $D$ be a circular brane diagram with any cut and let $f,g\in P_D$.  Then
for all $\lambda\in K_D$ and $\lambda^\dual\in K_{D^\dual}$ we have
\begin{equation}\label{eq:GM}
 \bigl\langle\eta_D(\lambda^\dual),\ \delta^D_{f,g}(\lambda)\bigr\rangle
 +\bigl\langle\eta_{D^\dual}(\lambda),\ \delta^{D^\dual}_{f^\dual,g^\dual}(\lambda^\dual)\bigr\rangle
 =0 .
\end{equation}
Equivalently, in $\operatorname{Hom}(K_D,E_D)=E_D\otimes Z_D$ we have
\begin{equation}\label{eq:GM-tensor}
 \delta^D_{f,g}+\tau\,(v_D\otimes u_D)\,\delta^{D^\dual}_{f^\dual,g^\dual}=0 ,
\end{equation}
where $\delta^{D^\dual}_{f^\dual,g^\dual}\in E_{D^\dual}\otimes Z_{D^\dual}$
and $\tau$ swaps the tensor factors.
\end{theorem}

The equivalence of \eqref{eq:GM} and \eqref{eq:GM-tensor} is a matter of
unwinding Definition~\ref{def:eta}. Namely, suppose we evaluate \eqref{eq:GM-tensor} on $(\sigma,\lambda)\in E_D^*\times K_D$ with $\sigma=\eta_D(\lambda^\dual)$. Then using $\pair{\sigma}{u_D(\beta^\dual_i)}=\lambda^\dual_i$ and $v_D^*=\eta_{D^\dual}$ we obtain \eqref{eq:GM}.
For the other direction recall that $\eta_D$ is an isomorphism, and hence \eqref{eq:GM} for all $\lambda^\dual$ implies~\eqref{eq:GM-tensor}.

Theorem \ref{thm:A} is proved in Sections~\ref{subsec:reduction}--\ref{subsec:proofA}.
First we look at the smallest interesting example.

\subsection{The Hilbert scheme of points of \texorpdfstring{$\C^2$}{C2}}
\label{subsec:Hilb}
Let $D$ be the circular diagram with one D5 brane $a$, one NS5 brane $z$, and
multiplicity $d$ on both arcs. Clearly, we have $D^\dual=D$. We present $D$ in normal form $z\,a$ with the cut between $a$ and $z$.

The combinatorial description of fixed points on $\X(D)$ are tie diagrams on $D$, and their associated signed indicator functions 
$\nu(-)=\nu(\cdot,1,1)$. It is easy to show that they are in bijection with partitions of $d$, using the following rules (see also Figure~\ref{fig:31tiediagram}). For a partition $\mu=(\mu_1\geq \mu_2 \geq \ldots)$, define $B(\mu)=\{\mu_r-r : r\geq 1\} \subset \Z$ (with $B(\emptyset)=\Z_{<0}$), and
\begin{equation}\label{eq:nu-mu}
 \nu_\mu(p)=\mathbf1_{p-1\in B(\emptyset)}-\mathbf1_{p-1\in B(\mu)}
 =\mathbf1_{p\le0}-\mathbf1_{p-1\in B(\mu)}.
\end{equation}

Define the \emph{boxes} of $\mu$ to be the pairs $(r,s)$ with $1\le s\le\mu_r$, and their {\em content} $c(r,s)=s-r$.
The {\em total content} of $\mu$ is 
\[
 T(\mu)=\sum_{\square\in\mu}c(\square)
 =\sum_{r\ge1}\sum_{s=1}^{\mu_r}(s-r)
 =\sum_{r\ge1}\Bigl[\binom{\mu_r+1}2-r\mu_r\Bigr].
\]
We leave it as a small exercise to prove that, in terms of $\nu$, the total content is the finite sum
\begin{equation}\label{eq:hilb_content}
T(\mu)=\dfrac12\sum_p p(1-p)\, \nu_\mu(p).
\end{equation}

\begin{lemma}\label{lem:hilb-weight}
For every $f\in P_D$,
\begin{equation}\label{eq:hilb-weight-ties}
 \ol{\wt(\ell|_f)}
 =\sum_p\nu_f(p)\Bigl[(1-p)a+\frac b4\,p(1-p)\Bigr]
 =da+\frac b4\sum_p p(1-p)\,\nu_f(p).
\end{equation}
Consequently, for partitions $\mu,\mu'$ of $d$ we have
\begin{equation}\label{eq:hilb-weight}
 \ol{\wt(\ell|_\mu)}=da+\frac b2\,T(\mu),
 \qquad
 \delta^D_{\mu,\mu'}(\ell)=\frac b2\bigl(T(\mu)-T(\mu')\bigr).
\end{equation}
\end{lemma}
\begin{proof}
By definition, $\wt(\ell|_f)$ is the sum of the weights of the tie portions
covering the lift of $I$ (Section~\ref{subsec:decoration}). By
Lemma~\ref{lem:intrinsic} the portion over $I^{(s)}$ has weight $a+\frac s2b$
modulo $\h$. Therefore, a left tie at $p$
contributes $(1-p)a+\frac b4p(1-p)$ and a right tie at $p$ contributes
$(p-1)a-\frac b4p(1-p)$; the sign of $\nu_f$ absorbs the difference, which
gives the first equality in \eqref{eq:hilb-weight-ties}.  

The second equality
is $\sum_p(1-p)\nu_f(p)=d$ from \eqref{eq:invariants}.  Now
\eqref{eq:hilb-weight} follows from the \eqref{eq:hilb_content} and the
definition \eqref{eq:delta} of $\delta$, notice the term $da$ cancels in the
difference, as it must since $\delta$ lands in $E_D=\Q b$.
\end{proof}

\medskip\noindent\emph{The mirror.}
We have $D^\dual=D$. Recalling 
\[
\nu_\mu(p)=\mathbf 1_{p\le0}-\mathbf 1_{p-1\in B(\mu)},
\qquad B(\mu)=\{\mu_u-u:u\ge1\}.
\] 
implies $B(\mu^t)=\{-1-b:b\notin B(\mu)\}$. This  gives $\nu_{\mu^t}(1-p)=-\nu_\mu(p)$, which match with the mirror bijection $\nu_{f^\dual}(1-p)=-\nu_f(p)$. Therefore, the tie mirror bijection is equivalent to transposing the partition. Since $p(1-p)$ is invariant under $p\mapsto1-p$, \eqref{eq:hilb_content} gives
\begin{equation}\label{eq:T-transpose}
 T(\mu^t)=-T(\mu).
\end{equation}

\medskip\noindent\emph{The GM pairings and the identity.}
The mirror map \eqref{eq:mirror-map} reads $b^\dual\mapsto-2\kappa$,
$\kappa^\dual\mapsto-\frac b2$ modulo $\h$, so by Lemma~\ref{lem:blocks}
and \eqref{eq:hilb-weight}
applied to $D$ and to $D^\dual$ we have
\[
 \bigl\langle\eta_D(\ell^\dual),\ \delta^D_{\mu,\mu'}(\ell)\bigr\rangle
 +\bigl\langle\eta_{D^\dual}(\ell),\ \delta^{D^\dual}_{\mu^t,{\mu'}^t}(\ell^\dual)\bigr\rangle
 =\bigl(T(\mu)-T(\mu')\bigr)+\bigl(T(\mu^t)-T({\mu'}^t)\bigr)=0 ,
\]
So for the Hilbert scheme, GM duality is just the
statement that transposition negates the total content of a partition.

\begin{example}\label{ex:31} \rm (C.f. Figure~\ref{fig:31tiediagram}.)
For $\mu=(3,1)$ one has $B(\mu)=\{2,-1,-3,-4,\dots\}$, so the only nonzero
values of $\nu_\mu$ are $\nu_\mu(-1)=1$ and $\nu_\mu(3)=-1$. Therefore $\sum\nu_\mu=0$ and
$\sum(1-p)\nu_\mu(p)=2+2=4=d$. We have 
$T(\mu)=\frac12\bigl[(-1)(2)-(3)(-2)\bigr]=2$, the sum of the contents
$0,1,2,-1$ of the four boxes. 

The corresponding tie diagram
has two ties at $a^{(0)}$, a left tie to $z^{(-1)}$ and a right tie to
$z^{(3)}$. In the proof of Lemma~\ref{lem:hilb-weight}
the left tie carries the weights $a-\frac b2,\,a$ and the right tie the
weights $a+\frac b2,\,a+b$, so $\ol{\wt(\ell|_\mu)}=4a+b$.  The transpose
$\mu^t=(2,1,1)$ has $B(\mu^t)=\{1,-1,-2,-4,-5,\dots\}$, hence
$\nu_{\mu^t}(-2)=1$, $\nu_{\mu^t}(2)=-1$, which is the reflection
$p\mapsto1-p$ of the previous indicator with the sign reversed, and
$T(\mu^t)=-2$.  
\end{example}
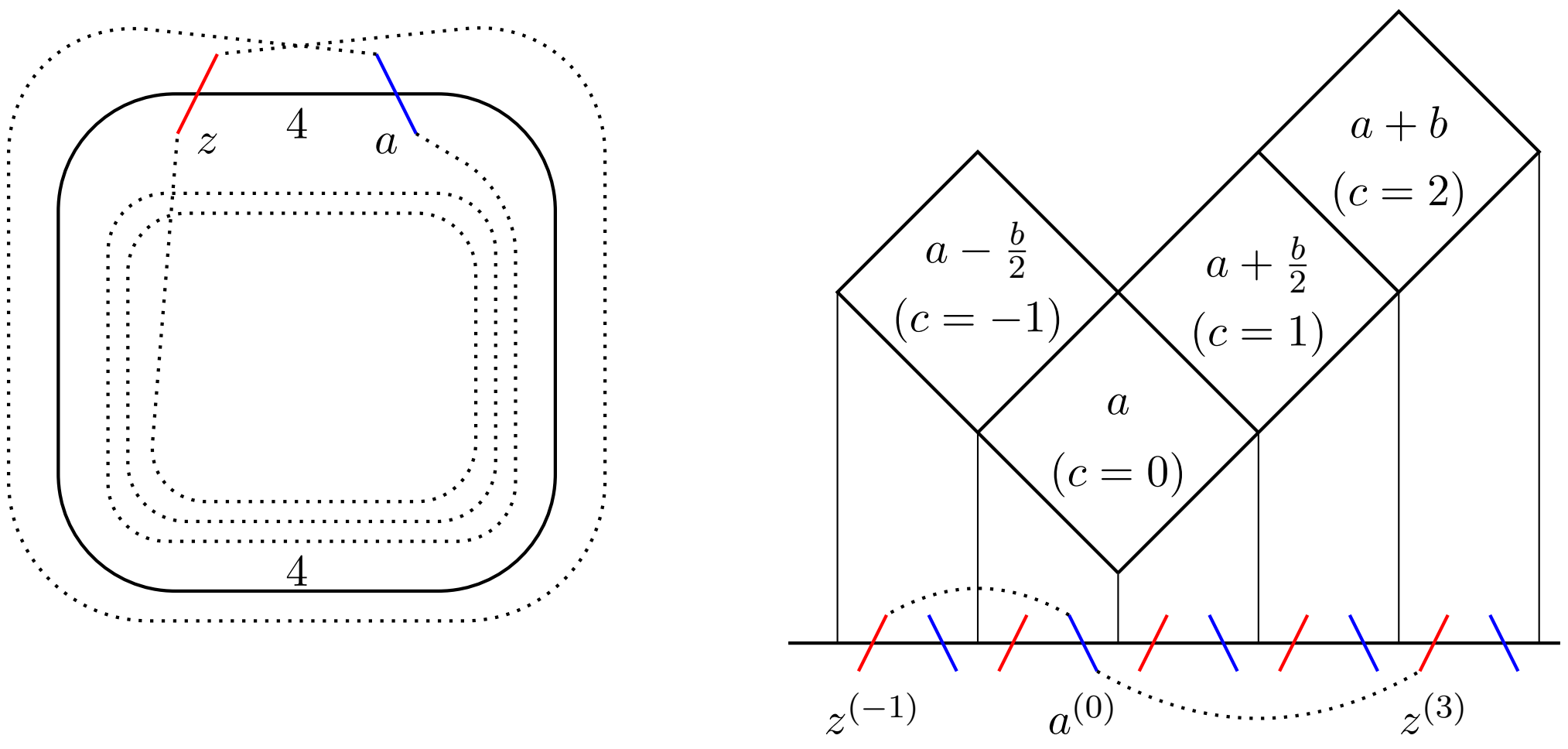
\begin{figure}[ht]
\[
\begin{tikzpicture}[scale=0.85]
    \draw[rounded corners=10mm, thick] (0,5) -- (5,5) -- (5,10) --(0,10) -- cycle;
    \draw[thick,red] (1.2,9.6)--(1.6,10.4);
    \draw[thick,blue] (3.6,9.6)--(3.2,10.4);
\draw[thick,dotted,rounded corners=5mm]
  (3.6,9.6) -- (4.6,9) -- (4.6,5.5) -- (.5,5.5) -- (.5,9) --
               (4.4,9) -- (4.4,5.7) -- (.7,5.7) -- (.7,8.8) --
               (4.2,8.8) -- (4.2,5.9) -- (.9, 5.9) -- (1.2,9.6);
\draw[thick,dotted,rounded corners=12mm]
  (3.2,10.4)-- (-.5,10.8) -- (-0.5,4.7) -- (5.5,4.7) -- (5.5,10.8) -- (1.6,10.4);
\draw (3.3,9.5) node {$a$};
\draw (1.5,9.5) node {$z$};
\draw (2.4,9.7) node {$4$}; \draw (2.4,5.2) node {$4$};
\end{tikzpicture}
\qquad\qquad
\begin{tikzpicture}[scale=.6,baseline=5]
    \draw[thick]  (-1,0) -- (10,0);
          \draw[thick,red] (0,-.4)--(0.4,0.4);
        \draw[thick,blue] (1.4,-.4)--(1,0.4);
          \draw[thick,red] (2,-.4)--(2.4,0.4);
        \draw[thick,blue] (3.4,-.4)--(3,0.4);
          \draw[thick,red] (4,-.4)--(4.4,0.4);
        \draw[thick,blue] (5.4,-.4)--(5,0.4);
          \draw[thick,red] (6,-.4)--(6.4,0.4);
        \draw[thick,blue] (7.4,-.4)--(7,0.4);
          \draw[thick,red] (8,-.4)--(8.4,0.4);
        \draw[thick,blue] (9.4,-.4)--(9,0.4);
\draw[thick,dotted] (3,.4) to[out=150,in=30] (0.4,.4);
\draw[thick,dotted] (3.4,-.4) to[out=-30,in=-150] (8,-.4);
\draw[thick, black] (3.7, 1) -- (-0.3,5) -- (1.7,7) -- (5.7,3);
\draw[thick, black] (3.7, 1) -- (9.7,7) -- (7.7,9) -- (1.7,3);
\draw[thick, black] (7.7,5) -- (5.7,7); 
\draw (3.2,-1) node {$a^{(0)}$};
\draw (0.2,-1) node {$z^{(-1)}$};
\draw (8.2,-1) node {$z^{(3)}$};
\draw (3.7,3.4) node {$a$}; \draw (3.7,2.4) node {$(c=0)$};
\draw (5.7,5.4) node {$a+\frac{b}{2}$}; \draw (5.7,4.4) node {$(c=1)$};
\draw (7.7,7.4) node {$a+b$}; \draw (7.7,6.4) node {$(c=2)$};
\draw (1.7,5.6) node {$a-\frac{b}{2}$}; \draw (1.7,4.6) node {$(c=-1)$};
\draw[thin] (3.7,1) -- (3.7,0);
\draw[thin] (5.7,3) -- (5.7,0);
\draw[thin] (7.7,5) -- (7.7,0);
\draw[thin] (9.7,7) -- (9.7,0);
\draw[thin] (1.7,3) -- (1.7,0);
\draw[thin] (-.3,5) -- (-.3,0);
\end{tikzpicture}
\]
\caption{The fixed point of $\operatorname{Hilb}^{4}(\C^2)$ corresponding to  
$\mu=(3,1)$, its tie diagram (left) with its universal cover (right). The box weights modulo $\h$ are shown.}
\label{fig:31tiediagram}
\end{figure}

\subsection{An example with four fixed points}
Take $D=z_1z_2a_1a_2$ with all four interval multiplicities equal to $1$.
We have four fixed points $f_{ij}$ ($i,j=1,2$) 
with $\nu(0,i,j)=1$, $\nu(1,i,j)=-1$ and all other entries zero.
The two determinant weights modulo $\h$ are
\[
\begin{array}{c|cc}
 f&\ell_1&\ell_2\\\hline
 f_{11}&a_1-b/2&a_1\\
 f_{12}&a_1+b/2&a_1\\
 f_{21}&a_2-b/2&a_2\\
 f_{22}&a_2+b/2&a_2.
\end{array}
\]
In the mirror normal form $f_{ij}^\dual=f_{ji}$.
Write $T_\lambda=\lambda_1+\lambda_2$ and
$T^\dual_\lambda=\lambda^\dual_1+\lambda^\dual_2$.  The parameter pairing gives
$\langle\eta_D(\lambda^\dual),a_1-a_2\rangle=\lambda^\dual_1$ and
$\langle\eta_D(\lambda^\dual),b\rangle=T^\dual_\lambda$, and likewise on the mirror.
For the pairs $(f_{11},f_{21})$ and $(f_{11},f_{12})$, respectively,
the two terms of \eqref{eq:GM} are
\[
 T_\lambda\lambda^\dual_1-\lambda^\dual_1T_\lambda=0,
 \qquad -\lambda_1T^\dual_\lambda+T^\dual_\lambda\lambda_1=0.
\]
All other differences follow by addition.  

\subsection{Reduction to normal form}
\label{subsec:reduction}
The first observation is that both sides of \eqref{eq:GM} are independent of the presentation of $D$.

\begin{lemma}\label{lem:reduction}
It suffices to prove Theorem~\ref{thm:A} for $D$ in normal form with the cut
between $a_n$ and~$z_1$.
\end{lemma}
\begin{proof}
If $r$ D5 and $s$ NS5 branes cross the cut, and
$\chi_i,\eta_j$ are the corresponding crossing indicators, recall the coordinate
changes\eqref{eq:cut-counts}
\[
 a_i'=a_i+m\chi_i t_1,\quad R_i'=R_i+s-m\chi_i,
 \qquad z_j'=z_j+\eta_j\zeta,\quad L_j'=L_j-r+n\eta_j.
\]
Thus $A_i'=A_i+s t_1$ and
$z_j'-L_j'\kappa=z_j-L_j\kappa+r\kappa$. Notice $A_i'-A_j'=A_i-A_j$ remains unchanged, so does the $\beta$. Recall \eqref{lem:eta}, the GM identity is cut independent.

\emph{Hanany--Witten transitions.}  Let $D_1,D_2$ differ by one transition,
with the D5 brane $a$ of $D_1$ becoming $\widetilde a=a+t_1$ in $D_2$, and let
$f\mapsto\Psi(f)$ be the induced bijection of fixed points.  

Since the
isomorphism $\X(D_1)\cong\X(D_2)$ is equivariant after this reparametrization
and matches the tautological bundles as in \eqref{eq:HW-taut}, 

An HW transition moving a D5 brane $a$ to the right
across an NS5 brane $z$. Write $\widetilde a=a+t_1$, let
$\Psi$ be the induced fixed-points bijection, and put
$c=d_1-d_2=\operatorname{lch}(a)$. 

Given $(\lambda_1,...,\lambda_n)\in K_D$, we are going to compare \[wt(\sum \lambda_il_i)|_f-wt(\sum \lambda_il_i)|_g\] and \[wt(\sum \lambda_i\tilde{l}_i)|_f-wt(\sum \lambda_i\tilde{l}_i)|_g\] for the determinant generators $l_i$ and $\tilde{l}_i$ attach to $\beta_i$ before and after HW transition.

To do so we only need to compare the determinant generator attached to $z$, which is represented by
$\det\xi_3$ before the move and by $\det\widetilde\xi_2$
afterwards. Taking determinants in ~\ref{eq:HW-taut} and applying ~\ref{lem:jump} gives
\begin{align*}
wt(det\tilde{\xi}_2)|_{\Psi(f)}-wt(det\xi_3)|_f
 &=wt(det\xi_1)|_f-wt(det\xi_2)|_f+\widetilde a-c t_2\\
 &=c a-\binom{1-c}{2}\hbar+a+t_1-c t_2.
\end{align*}
Which is fixed points independent.

Also, $\eta_D(\lambda^!)$ and $\eta_{D^!}(\lambda)$ are HW invariant since $A$ is HW invariant and the cocharacter $\sigma=\eta_D(\lambda^!)$ is defined such that ~\ref{lem:eta} holds. Hence we conclude that ~\ref{thm:A} is HW invariance.

\end{proof}

For the rest of this section $D$ is in normal form with the cut between $a_n$
and $z_1$, with the lifts and cells of Section~\ref{subsec:cells}.  Then
$R_i=L_j=0$, $A_i=a_i$, the lifted K\"ahler coordinates are
$\widetilde z_j^{(s)}=z_j+s\zeta$, and $\beta_j$ is as in \eqref{eq:beta}.
Write $S=\sum_jz_j$ and $S_{>j}=\sum_{j'>j}z_{j'}$.

\subsection{Linearity in the tie diagram}
\label{subsec:linearity}
The relative restriction maps $\delta^D_{f,g}$ are linear in $\nu_f$, with an explicit
coefficient for each cell.  It is convenient to package
$\delta^D_{f,g}$ through the $Z_D$-valued combination
\begin{equation}\label{eq:Wf}
 U_D(f):=2\sum_{j=1}^m\beta_j\,\ol{\wt(\ell_j|_f)}\ \in\ Z_D\otimes X^*(\A_D)_\Q ,
 \quad
 U_D(f)-U_D(g)=2\sum_j\beta_j\otimes\delta^D_{f,g}(\ell_j)=2\,\delta^D_{f,g}
\end{equation}
under $Z_D\otimes E_D=\operatorname{Hom}(K_D,E_D)$.  

\begin{lemma}\label{lem:Psi}
For $(p,i,j)\in\Z\times\{1,\dots,n\}\times\{1,\dots,m\}$, put
\begin{equation}\label{eq:Psi}
 \Psi_D(p,i,j)=2a_i\bigl(z_j-z_1+(p-1)\zeta\bigr)
 -b\Bigl[\bigl(m(1-p)-j\bigr)z_j-S_{>j}+pS-\tfrac m2\,p(p-1)\,\zeta\Bigr].
\end{equation}
Then for every $f\in P_D$, the following equation holds:
\begin{equation}\label{eq:U-linear}
 U_D(f)=\sum_{p,i,j}\nu_f(p,i,j)\,\Psi_D(p,i,j).
\end{equation}
\end{lemma}

\begin{proof}
By Lemma~\ref{lem:intrinsic}, for a given tie, over the lifted intervals
$I=I_{j'}^{(s')}$ it covers, the tie at $(p,i,j)$ contributes to $U_D(f)$
the sum $2\sum_I\beta(I)\,\ol{x_t(I)}$, where $\beta(I)=\beta_{j'}$ and $x_t(I)$ is the
weight of the portion over $I$.

\emph{Left ties}, $p\le0$.  The tie runs from $z_j^{(p)}$ to $a_i^{(0)}$.  It
covers $I_{j'}^{(p)}$ for $j'\ge j$ and $I_{j'}^{(s')}$ for all $j'$ and
$p<s'\le0$.  The number of NS5 lifts between $I_{j'}^{(s')}$ and $a_i^{(0)}$
is $k=m(1-s')-j'$.  So the contribution is $2Xa_i-bY$ with,
\begin{align*}
 X&=\sum_{j'\ge j}\beta_{j'}+\sum_{s'=p+1}^{0}\sum_{j'}\beta_{j'}
  =\widetilde z_j^{(p)}-\widetilde z_1^{(1)}=z_j-z_1+(p-1)\zeta,\\
 Y&=\sum_{j'\ge j}\beta_{j'}\bigl(m(1-p)-j'\bigr)
   +\sum_{s'=p+1}^{0}\sum_{j'}\beta_{j'}\bigl(m(1-s')-j'\bigr).
\end{align*}
With $\sum_{j'\ge j}\beta_{j'}=z_j-z_1-\zeta$,
$\sum_{j'\ge j}j'\beta_{j'}=jz_j+S_{>j}-mz_1-m\zeta$,
$\sum_{j'}\beta_{j'}=-\zeta$, $\sum_{j'}j'\beta_{j'}=S-mz_1-m\zeta$,
$\sum_{s'=p+1}^0(1-s')=\frac12p(p-1)$ and $\sum_{s'=p+1}^01=-p$,
\[
 Y=m(1-p)(z_j-z_1-\zeta)-(jz_j+S_{>j}-mz_1-m\zeta)-\tfrac m2p(p-1)\zeta
   +p(S-mz_1-m\zeta),
\]
in which the coefficients of $z_1$ and of $\zeta$ vanish
($-m(1-p)+m-mp=0$), leaving
$Y=(m(1-p)-j)z_j-S_{>j}+pS-\frac m2p(p-1)\zeta$.  So a left tie contributes
$\Psi_D(p,i,j)$.

\emph{Right ties}, $p\ge1$.  The tie runs from $a_i^{(0)}$ to $z_j^{(p)}$.  It
covers $I_{j'}^{(s')}$ for all $j'$ and $1\le s'\le p-1$, and $I_{j'}^{(p)}$
for $j'<j$.  The number of NS5 lifts between $a_i^{(0)}$ and $I_{j'}^{(s')}$ is
$k=m(s'-1)+j'$, and the $b$-term has the opposite sign.  The contribution is
$2X'a_i+bY'$ with,
\[
 X'=\sum_{s'=1}^{p-1}\sum_{j'}\beta_{j'}+\sum_{j'<j}\beta_{j'}
   =\widetilde z_1^{(1)}-\widetilde z_j^{(p)}=-X,
\]
\begin{align*}
 Y'&=\sum_{s'=1}^{p-1}\sum_{j'}\beta_{j'}\bigl(m(s'-1)+j'\bigr)
    +\sum_{j'<j}\beta_{j'}\bigl(m(p-1)+j'\bigr)\\
   &=-\tfrac m2(p-1)(p-2)\zeta+(p-1)(S-mz_1-m\zeta)
     +m(p-1)\bigl(-\zeta-(z_j-z_1-\zeta)\bigr)\\
   &\qquad+\bigl(S-mz_1-m\zeta-(jz_j+S_{>j}-mz_1-m\zeta)\bigr)=Y.
\end{align*}
Here we use $\sum_{j'<j}\beta_{j'}=-\zeta-(z_j-z_1-\zeta)$,
$\sum_{j'<j}j'\beta_{j'}=(S-mz_1-m\zeta)-(jz_j+S_{>j}-mz_1-m\zeta)$ and
$\sum_{s'=1}^{p-1}(s'-1)=\frac12(p-1)(p-2)$.  So a right tie contributes
$-2Xa_i+bY=-\Psi_D(p,i,j)$, the same polynomial with the sign of $\nu_f$.
\end{proof}

\subsection{The mirror side and the cell coefficient}
Present $D^\dual$ in normal form with the lifts in Lemma~\ref{lem:mirror}. Its cells are $(p^\dual,j,i)$ with $\nu_{f^\dual}(1-p,j,i)=-\nu_f(p,i,j)$.
By Lemma~\ref{lem:Psi} applied to $D^\dual$ (with $(n,m)$ replaced by
$(m,n)$) we have
\[
 U_{D^\dual}(f^\dual)=\sum_{p,i,j}\nu_{f^\dual}(1-p,j,i)\,\Psi_{D^\dual}(1-p,j,i)
 =-\sum_{p,i,j}\nu_f(p,i,j)\,\Psi_{D^\dual}(1-p,j,i).
\]
Modulo $\h$, since $L_j=R_i=0$ when $D$ is normal, the mirror map reads
$a^\dual_j\mapsto z_j$, $z^\dual_i\mapsto a_i$, $b^\dual\mapsto-2\kappa$,
$\kappa^\dual\mapsto-\frac12b$, $\zeta^\dual=m\kappa^\dual\mapsto-\frac m2b$.

\begin{remark}[Deck shifts]\label{rem:deck} \rm
The lifts $a^{\dual(0)}_j$ of Lemma~\ref{lem:mirror} sit one period to the
right of the lifts used in Definition~\ref{def:mirror-fixed-point}, so strictly
speaking $\varphi_D$ sends the coordinate $a^\dual_j$ of these lifts to
$z_j$ plus a common deck shift in $\Q\zeta\subset Z_D$.  Such a common shift
changes $\varphi_D(\Psi_{D^\dual}(1-p,j,i))$ by
$2\zeta\,\varphi_D\bigl(z^\dual_i-z^\dual_1+(p^\dual-1)\zeta^\dual\bigr)
=2\zeta(a_i-a_1)+nm\,p\,\kappa b$, which is of the form
$\varrho'_i+\gamma'p$, and is indeed harmless in the argument below (see
the proof of Theorem~\ref{thm:A}).  We suppress it.
\end{remark}

\begin{proposition}\label{prop:Lambda}
Let
$\Lambda(p,i,j)=\Psi_D(p,i,j)-\varphi_D\bigl(\Psi_{D^\dual}(1-p,j,i)\bigr)$,
reduced modulo $\h$.  Then
\begin{equation}\label{eq:Lambda}
 \Lambda(p,i,j)=\varrho_i+\varsigma_j+\gamma\,p+c ,
\end{equation}
with
\begin{align*}
 \varrho_i&=-2z_1a_i+2(i-n)\kappa a_i+2\kappa\sum_{i'>i}a_{i'},
 &\varsigma_j&=2a_1z_j+(j-m)bz_j+b\sum_{j'>j}z_{j'},\\
 \gamma&=2\kappa\sum_{i'}a_{i'}-b\sum_{j'}z_{j'},
 &c&=-2\kappa\sum_{i'}a_{i'}.
\end{align*}
\end{proposition}
\begin{proof}
With $p^\dual=1-p$, one has $p^\dual-1=-p$, $1-p^\dual=p$, and
$p^\dual(p^\dual-1)=p(p-1)$, so
\[
 \varphi_D\bigl(\Psi_{D^\dual}(1-p,j,i)\bigr)
 =2z_j\Bigl(a_i-a_1+\frac{pm}2b\Bigr)
 +2\kappa\Bigl[(np-i)a_i-\sum_{i'>i}a_{i'}+(1-p)\sum_{i'}a_{i'}
 +\frac{nm}4p(p-1)b\Bigr],
\]
while, with $\zeta=n\kappa$,
\[
 \Psi_D(p,i,j)=2a_iz_j-2z_1a_i+2(p-1)n\kappa a_i-\bigl(m(1-p)-j\bigr)bz_j
 +bS_{>j}-pbS+\frac{mn}2p(p-1)b\kappa .
\]
After subtracting $\varphi_D\bigl(\Psi_{D^\dual}(1-p,j,i)\bigr)$, $2a_iz_j$ cancels, the quadratic terms
$\frac{mn}2p(p-1)b\kappa$ cancel, the coefficient of $\kappa a_i$ is
$2(p-1)n-2(np-i)=2(i-n)$ and that of $bz_j$ is $-(m(1-p)-j)-pm=j-m$, both
independent of $p$. And the remaining terms gives \eqref{eq:Lambda}.
\end{proof}

\subsection{Proof of Theorem~\ref{thm:A}}
\label{subsec:proofA}
By Lemma~\ref{lem:reduction}, we can assume $D$ in normal form.  The
difference $U_{D^\dual}(f^\dual)-U_{D^\dual}(g^\dual)$ lies in
$Z_{D^\dual}\otimes E_{D^\dual}$. And by Lemma~\ref{lem:blocks}, it is easy to see the
mirror map $\varphi_D$, reduced modulo $\h$, descend to a map:
\[\tau(u_D\otimes v_D):Z_{D^\dual}\otimes E_{D^\dual}\to Z_D\otimes E_D.\]  Hence, by
\eqref{eq:Wf} on both sides and Lemma~\ref{lem:Psi} on both sides,
\begin{align*}
 2\Bigl(\delta^D_{f,g}+\tau(v_D\otimes u_D)\delta^{D^\dual}_{f^\dual,g^\dual}\Bigr)
 &=\bigl(U_D(f)-U_D(g)\bigr)+\varphi_D\bigl(U_{D^\dual}(f^\dual)-U_{D^\dual}(g^\dual)\bigr)\\
 &=\sum_{p,i,j}(\nu_f-\nu_g)(p,i,j)\Bigl[\Psi_D(p,i,j)-\varphi_D\bigl(\Psi_{D^\dual}(1-p,j,i)\bigr)\Bigr]\\
 &=\sum_{p,i,j}(\nu_f-\nu_g)(p,i,j)\,\Lambda(p,i,j),
\end{align*}
where the minus sign in front of $\varphi_D(\Psi_{D^\dual})$ comes from
\eqref{eq:N-mirror}.  By Proposition~\ref{prop:Lambda} the right side is
\[
 \sum_i\varrho_i\sum_{p,j}(\nu_f-\nu_g)
 +\sum_j\varsigma_j\sum_{p,i}(\nu_f-\nu_g)
 +\gamma\sum_{p,i,j}p\,(\nu_f-\nu_g)
 +c\sum_{p,i,j}(\nu_f-\nu_g),
\]
and each of the four sums vanishes by \eqref{eq:conservation}.  This proves
\eqref{eq:GM-tensor}, hence \eqref{eq:GM}.  \qed

\subsection{The absolute form}
\label{subsec:absolute}
The proof of Theorem~\ref{thm:A} evaluates the functional
$\nu\mapsto\sum\nu\Lambda$ on differences $\nu_f-\nu_g$, on which the three
invariants \eqref{eq:invariants} cancel. And the
same functional produces a constant when evaluate on $\nu_f$ itself.  Recall that
$\sum_{p,j}\nu_f(p,i,j)=\lch(a_i)$, $\sum_{p,i}\nu_f(p,i,j)=\lch(z_j)$ and
$\sum_{p,i,j}(1-p)\nu_f(p,i,j)=d_{I_m}$. Define
\[
 A:=\sum_{i=1}^n\lch(a_i)=\sum_{j=1}^m\lch(z_j),
\]
the two sums are equal because the D3 multiplicity returns to its initial
value after one turn around the circle.

\begin{proposition}[Absolute form]\label{prop:absolute}
Let $D$ be in normal form with the cut between $a_n$ and $z_1$, in the
coordinates of Proposition~\ref{prop:Lambda}.  For every $f\in P_D$,
\begin{equation}\label{eq:Omega}
 U_D(f)+\varphi_D\bigl(U_{D^\dual}(f^\dual)\bigr)
 =\Omega_D:=\sum_{i=1}^n\lch(a_i)\,\varrho_i+\sum_{j=1}^m\lch(z_j)\,\varsigma_j
 +\gamma\,(A-d_{I_m})+c\,A 
\end{equation}
is an element of
$\operatorname{span}_\Q\{z_j,\kappa\}\otimes\operatorname{span}_\Q\{a_i,b\}$
which does not depend on $f$.  In particular, if all local charges of $D$
vanish, and all multiplicities equal $d:=d_{I_m}$, then,
\begin{equation}\label{eq:Omega-balanced}
 \Omega_D=-d\Bigl(2\kappa\otimes\sum_{i=1}^na_i-\sum_{j=1}^mz_j\otimes b\Bigr).
\end{equation}
\end{proposition}
\begin{proof}
By Lemma~\ref{lem:Psi} on both sides, \eqref{eq:N-mirror} and
Proposition~\ref{prop:Lambda},
\[
 U_D(f)+\varphi_D\bigl(U_{D^\dual}(f^\dual)\bigr)
 =\sum_{p,i,j}\nu_f(p,i,j)\,\Lambda(p,i,j)
 =\sum_i\varrho_i\sum_{p,j}\nu_f
 +\sum_j\varsigma_j\sum_{p,i}\nu_f
 +\gamma\sum_{p,i,j}p\,\nu_f
 +c\sum_{p,i,j}\nu_f .
\]
Now insert the values \eqref{eq:invariants}, the row sums are $\lch(a_i)$, the column sums
are $\lch(z_j)$, $\sum\nu_f=A$, and
$\sum p\,\nu_f=\sum\nu_f-\sum(1-p)\nu_f=A-d_{I_m}$.  This is
\eqref{eq:Omega}.  If all local charges vanish, then only
$\gamma(A-d)+cA=-d\gamma$ survives, hence
$\gamma=2\kappa\sum_ia_i-b\sum_jz_j$.
\end{proof}
Unlike
the relative identity, \eqref{eq:Omega} is sensitive to the common deck
shift of Remark~\ref{rem:deck}, which under a shift by $s\zeta$ adds to $\Omega_D$ the term
$-s\bigl(2\zeta\sum_i(a_i-a_1)\lch(a_i)+nm\,(A-d_{I_m})\,\kappa b\bigr)$.  This term is the image under $\varphi_D$ of a common shift
of the $a^\dual_j$, and it disappears under the normalization introduced next in the balanced case.

\emph{The traceless normalization.}  The weight $\ol{\wt(\lambda|_f)}$ has
$a$-degree $\sum_j\lambda_jd_{I_j}$ and so does not lie in $E_D$.  Only its
differences do.  The simplest way to force it into $E_D$ is to pass from
the maximal torus of $GL_n$ to that of $SL_n$.  For
$w=\sum_ic_ia_i+pt_1+qt_2\in X^*(\T_D)_\Q$ put $\deg_aw=\sum_ic_i$ and
\begin{equation}\label{eq:traceless}
 w^\circ=\ol w-\frac{\deg_aw}{n}\sum_{i=1}^na_i\ \in E_D .
\end{equation}
It defines a projection $\pi_a\colon X^*(\A_D)_\Q\to E_D$ with kernel $\Q\sum_ia_i$.
Similarily on the mirror side, switch $m$ and $n$, we have projection $\pi_{a^\dual}\colon X^*(\A_{D^\dual})_\Q\to E_{D^\dual}$.  Note that
$\delta^D_{f,g}(\lambda)=\wt(\lambda|_f)^\circ-\wt(\lambda|_g)^\circ$, because $w\mapsto w^\circ$ kills every common shift of the $a_i$.

\begin{corollary}[Absolute GM identity for balanced diagrams]\label{cor:absolute}
Let $D$ be in normal form and suppose that all local charges of $D$ vanish.
Then for all $f\in P_D$, $\lambda\in K_D$ and $\lambda^\dual\in K_{D^\dual}$,
\begin{equation}\label{eq:GM-absolute}
 \bigl\langle\eta_D(\lambda^\dual),\ \wt(\lambda|_f)^\circ\bigr\rangle
 +\bigl\langle\eta_{D^\dual}(\lambda),\ \wt(\lambda^\dual|_{f^\dual})^\circ\bigr\rangle
 =0 .
\end{equation}
\end{corollary}

For general $D$, the same computation expresses the failure of the absolute
identity in the traceless normalization by the single tensor
$(\pi_z\otimes\pi_a)(\Omega_D)\in Z_D\otimes E_D$, which need not vanish.  We
return to this in Section~\ref{subsec:GM-comparison}.

\section{GM duality on the tautological model}
\label{sec:GMduality}

\subsection{The statement}
Let
\[
 H_D=E_D\oplus K_D,\qquad H_{D^\dual}=E_{D^\dual}\oplus K_{D^\dual},
\]
and let $B\colon H_D\times H_{D^\dual}\to\Q$ be the bilinear form
\begin{equation}\label{eq:B}
 B\bigl((\varepsilon,\lambda),(\varepsilon^\dual,\lambda^\dual)\bigr)
 =\pair{\eta_D(\lambda^\dual)}{\varepsilon}+\pair{\eta_{D^\dual}(\lambda)}{\varepsilon^\dual}.
\end{equation}
Assume $P_D\ne\varnothing$ and fix a base point $f_0\in P_D$.  For $f\in P_D$, let
$\rho_f\colon H_D\to E_D$ $\rho_f(\varepsilon,\lambda)=\varepsilon+\delta^D_{f,f_0}(\lambda)$,
be the \emph{restriction to $f$ normalized by the equivariant weight at $f_0$}, i.e. we consider $\delta^D_{f,f_0}$ instead of $\delta^D_f$, and put
\begin{equation}\label{eq:Vf}
 V_f=\ker\rho_f=\gr\bigl(-\delta^D_{f,f_0}\bigr)
 =\{(-\delta^D_{f,f_0}(\lambda),\lambda):\lambda\in K_D\}\subset H_D,
\end{equation}
and similarly $V^\dual_{f^\dual}=\ker\rho^\dual_{f^\dual}\subset H_{D^\dual}$
with base point $f_0^\dual$.

\begin{theorem}[GM duality on the tautological model]\label{thm:C}
Assume $P_D\ne\varnothing$ and choose $f_0\in P_D$. Then the pairing $B$ is perfect, under such pairing we have $E_D\perp E_{D^\dual}$ and $K_D\perp K_{D^\dual}$,
and for every $f\in P_D$,
\[
 \bigl(V_f\bigr)^{\perp_B}=V^\dual_{f^\dual}.
\]
\end{theorem}
\begin{proof}
$B$ is block off-diagonal.  Its $E_D\times K_{D^\dual}$ block is $\eta_D$ and
its $K_D\times E_{D^\dual}$ block is $\eta_{D^\dual}$, both are isomorphisms. The other two blocks vanish, which implies
$E_D\perp E_{D^\dual}$, $K_D\perp K_{D^\dual}$.  Hence $B$ is perfect.  For $\lambda\in K_D$,
$\lambda^\dual\in K_{D^\dual}$,
\[
 B\bigl((-\delta^D_{f,f_0}\lambda,\lambda),(-\delta^{D^\dual}_{f^\dual,f_0^\dual}\lambda^\dual,\lambda^\dual)\bigr)
 =-\pair{\eta_D\lambda^\dual}{\delta^D_{f,f_0}\lambda}
  -\pair{\eta_{D^\dual}\lambda}{\delta^{D^\dual}_{f^\dual,f_0^\dual}\lambda^\dual}
 =0
\]
by Theorem~\ref{thm:A} with $g=f_0$.  So $V^\dual_{f^\dual}\subset V_f^\perp$,
and $\dim V_f^\perp=(n+m)-m=n=\dim V^\dual_{f^\dual}$.
\end{proof}

The computation above is reversible, hence we have,

\begin{lemma}[Converse]\label{lem:converse}
Let $B'$ be any perfect pairing $H_D\times H_{D^\dual}\to\Q$ with
$E_D\perp E_{D^\dual}$ and $K_D\perp K_{D^\dual}$, and let $\eta'_D:K_D\to E_{D^\dual}$,
$\eta'_{D^\dual}:K_{D^\dual}\to E_D$ be its two nonzero blocks.  Then $B'(V_f,V^\dual_{f^\dual})=0$
for all $f\in P_D$ if and only if \eqref{eq:GM} holds for
$\eta'_D,\eta'_{D^\dual}$.  In particular,
Theorem~\ref{thm:A} and Theorem~\ref{thm:C} are equivalent.
\end{lemma}

\begin{remark} \rm
The base point is just a convenience.
The subspaces $V_f$ depend on $f_0$, but
replacing $f_0$ by $f_1$ transforms $H_D$ by the shear
$S(\varepsilon,\lambda)=(\varepsilon-\delta^D_{f_0,f_1}(\lambda),\lambda)$. For $x=(\varepsilon,\lambda)$ and
$y=(\varepsilon^\dual,\lambda^\dual)$, Theorem~\ref{thm:A} gives,
\[
 B(Sx,S^\dual y)-B(x,y)
 =-\pair{\eta_D\lambda^\dual}{\delta^D_{f_0,f_1}\lambda}
 -\pair{\eta_{D^\dual}\lambda}{\delta^{D^\dual}_{f_0^\dual,f_1^\dual}\lambda^\dual}=0 .
\]
Thus $B$ is invariant under the change of base point, and this invariance is
another equivalent form of Theorem~\ref{thm:A}.  We return to this in
Section~\ref{subsec:GM-comparison}, where the base point normalization is replaced by
the fixed-point average of \cite{GM10}. And for balanced diagrams, another normalization called
traceless normalization of Section~\ref{subsec:absolute} is also applicable.
\end{remark}

\subsection{Comparison with Goresky--MacPherson duality}
\label{subsec:GM-comparison}
We recall the original formulation by Goresky and MacPherson and explain in what sense
Theorem~\ref{thm:C} is a version of it.  Let $A$ be a torus acting on a
smooth connected variety $\X$ with a nonempty finite fixed-points set, and write
$W=H^2_A(\X;\Q)$.  The inclusion $\mathfrak t^*=H^2_A(\mathrm{pt};\Q)\to W$ split by the restriction $r_p\colon W\to H^2_A(p;\Q)=\mathfrak t^*$ for any
fixed point $p$, and its kernel,
\[
 K_p=\ker\bigl(H^2_A(\X)\to H^2_A(p)\bigr)
\]
is a complement of $\mathfrak t^*$ in $W$.  The arrangement of the
$\dim\mathfrak t$-dimensional subspaces $\ann(K_p)\subset W^*=H_2^A(\X)$ is
the degree-two part of the picture of $\operatorname{Spec}H^*_A(\X)$ studied
in \cite{GM10}, and they are equal when $\operatorname{Spec}H^2_A(\X)$ generates $H^*_A(\X)$. 
\begin{definition}[GM duality]\label{def:GM}
Let $\X,\X^\dual$ be as above, with tori $A,A^\dual$ and a bijection
$p\mapsto p^\dual$ of fixed points.  A perfect pairing
$\langle\,,\rangle\colon H^2_A(\X)\times H^2_{A^\dual}(\X^\dual)\to\Q$
establishes \emph{Goresky--MacPherson duality} if,
\begin{enumerate}
\item the equivariant parameters of the two sides are orthogonal,
 $\mathfrak t^*\perp(\mathfrak t^\dual)^*$, and
\item for every fixed point $p$, $K_p^\perp=K_{p^\dual}$.
\end{enumerate}
\end{definition}
Conditions (1) and (2) together imply that the annihilator of
$\mathfrak t^*$ is $(\mathfrak t^\dual)^*$.  Indeed, by dimensions, condition
(2) gives $\dim W=\dim\mathfrak t+\dim\mathfrak t^\dual$.
Consequently $\mathfrak t^*$ pairs perfectly with
$H^2_{A^\dual}(\X^\dual)/(\mathfrak t^\dual)^*\cong H^2(\X^\dual)$.  The
equivariant para\-meters of $\X$ are the K\"ahler parameters of $\X^\dual$.
Condition (2) is the duality of the two arrangements observed in
\cite{GM10} for type $A$ flag varieties and the corresponding slices. Both conditions are listed among the expected properties of
symplectic dual pairs in \cite[\S10]{BLPW16}.
\begin{remark}
When identifying the quotient by equivariant parameters with
ordinary degree-two cohomology, we assume on both sides the exact sequence
$0\to H_A^2(\mathrm{pt})\to H_A^2(\X)\to H^2(\X)\to0$ holds.
The arrangement definition itself does not require this extra hypothesis.
\end{remark}
Now let $\X=\X(D)$, $A=\A_D/(\text{diagonal})$ the symplectic torus modulo the framing diagonal, so that $\mathfrak t^*=E_D$, and let $f_0\in P_D$.  Each determinant
line bundle $\det\xi_{I_j}$ is $\T_D$-equivariant, twisting by $\det\xi_{I_j}|^{-1}_{f_0}$ makes the diagonal act trivially, and the twisted bundle is $A$-equivariant, then we obtain a linear map
\begin{equation}\label{eq:Phi}
 \Phi_{f_0}\colon H_D=E_D\oplus K_D\longrightarrow H^2_A(\X(D);\Q),
 \qquad
 (\varepsilon,\lambda)\mapsto\varepsilon+\sum_j\lambda_j\,
 c_1^A\bigl(\det\xi_{I_j}\otimes\C_{-\wt(\det\xi_{I_j}|_{f_0})}\bigr).
\end{equation}
By construction $r_f\circ\Phi_{f_0}=\rho_f$ for every $f\in P_D$. And we have
$\Phi_{f_0}(V_f)\subset K_f$.

\begin{proposition}\label{prop:GM-genuine}
Suppose that $\Phi_{f_0}$ and $\Phi^\dual_{f_0^\dual}$ are isomorphisms.
Then the pairing on $H^2_A(\X(D))\times H^2_{A^\dual}(\X(D^\dual))$
transported from $B$ establishes GM duality for the pair
$(\X(D),\X(D^\dual))$, and it does
not depend on the base points.
\end{proposition}
\begin{proof}
If $\Phi_{f_0}$ is an isomorphism, then $\Phi_{f_0}(V_f)=K_f$ for all $f$ and
$\Phi_{f_0}(E_D)=\mathfrak t^*$. Then it is easy to see that under the transported pairing
 $t^*\perp (t^!)^*$ and
$   K_f^\perp=K^\dual_{f^\dual}$. Changing $f_0$ to
$f_1$ replaces $\Phi_{f_0}$ by $\Phi_{f_1}=\Phi_{f_0}\circ S^{-1}$ with the
shear $S$, since $B$ is invariant under such shear.  So the
transported pairing is independent of the base points.
\end{proof}

\medskip
\emph{Base-point-free forms.} In \cite{GM10} they do not choose any particular base point.  
They use the \emph{average}
$\mu_0=\frac1{|\X^A|}\sum_p\mu_p$ of
$\mu_p\colon\mathfrak t=H_2^A(p)\to H_2^A(\X)$ defined by the fixed points
\cite[(3.2)]{GM10} to split $0\to H_2(\X)\to H_2^A(\X)\to\mathfrak t\to0$.  

There are two natural base-point-free normalizations of the tautological model.  The first is the dual of the average of
\cite{GM10},
\begin{equation}\label{eq:Phi-av}
 \Phi_{\mathrm{av}}=\frac1{|P_D|}\sum_{f_0\in P_D}\Phi_{f_0},
 \qquad
 \Phi_{\mathrm{av}}(\varepsilon,\lambda)=\varepsilon+\sum_j\lambda_j\,
 c_1^A\Bigl(\det\xi_{I_j}\otimes\C_{-\frac1{|P_D|}\sum_{g\in P_D}\wt(\det\xi_{I_j}|_g)}\Bigr),
\end{equation}
each determinant bundle being twisted by minus of its fixed-point weights average.  The corresponding subspaces are
$V^{\mathrm{av}}_f=\gr(-\delta^D_{f,\mathrm{av}})$ with
$\delta^D_{f,\mathrm{av}}=\frac1{|P_D|}\sum_{g\in P_D}\delta^D_{f,g}$.

The second is the traceless normalization \eqref{eq:traceless} of
Section~\ref{subsec:absolute}. After restricting our torus to the subtorus $\{\prod_ia_i=1\}$ of the framing
torus, $\Phi$ and $V_f$ reduces to,
\begin{equation}\label{eq:Phi-tr}
 \Phi^\circ(\varepsilon,\lambda)=\varepsilon+\sum_j\lambda_j\,
 c_1^A\Bigl(\det\xi_{I_j}\otimes\C_{-\frac{d_{I_j}}n\sum_ia_i}\Bigr),
 \qquad
 V^\circ_f=\{(-\wt(\lambda|_f)^\circ,\lambda):\lambda\in K_D\}.
\end{equation}
Both two normalizations are of the form $\Phi_{f_0}\circ S^{-1}$ for a shear
$S(\varepsilon,\lambda)=(\varepsilon-s(\lambda),\lambda)$ with
$s\in\Hom(K_D,E_D)$, more precisely $s=\delta^D_{f_0,\mathrm{av}}$ for
$\Phi_{\mathrm{av}}$ and $s=\wt(\,\cdot\,|_{f_0})^\circ$ for $\Phi^\circ$,
and in each case the new subspaces are $SV_f$.  On the mirror side
$S^\dual$ is defined by the corresponding $s^\dual\in\Hom(K_{D^\dual},E_{D^\dual})$.

We want to see whether the pairing $B$ is invariant under either of these two pairs of shears. As in the remark after
Lemma~\ref{lem:converse},
\begin{equation}\label{eq:shear-B}
 B(Sx,S^\dual y)-B(x,y)
 =-\pair{\eta_D(\lambda^\dual)}{s(\lambda)}-\pair{\eta_{D^\dual}(\lambda)}{s^\dual(\lambda^\dual)},
 \qquad x=(\varepsilon,\lambda),\ y=(\varepsilon^\dual,\lambda^\dual).
\end{equation}

\begin{proposition}[Base-point-free forms]\label{prop:basepoint-free}
Assume $P_D\ne\varnothing$.
\begin{enumerate}
\item[(1)] For every $D$, Theorem~\ref{thm:C} holds with the averaged
 subspaces, that is, $(a):(V^{\mathrm{av}}_f)^{\perp_B}=V^{\dual\,\mathrm{av}}_{f^\dual}$ for
 all $f\in P_D$.  If $\Phi_{f_0}$, $\Phi^\dual_{f_0^\dual}$ are
 isomorphisms, $(b):$ the pairing on $H^2_A(\X(D))\times H^2_{A^\dual}(\X(D^\dual))$ transported by $(\Phi_{\mathrm{av}},\Phi^\dual_{\mathrm{av}})$ is that of
 Proposition~\ref{prop:GM-genuine}.
\item[(2)] The same two statements $(a)$ and $(b)$ hold for the traceless subspaces
 $V^\circ_f$, $V^{\dual\circ}_{f^\dual}$ and the maps $\Phi^\circ$,
 $\Phi^{\dual\circ}$ if and only if the absolute identity
 \eqref{eq:GM-absolute} holds for one, equivalently for every $f\in P_D$.
 For $D$ in normal form, this is the condition
 $(\pi_z\otimes\pi_a)(\Omega_D)=0$ of Section~\ref{subsec:absolute}.  In
 particular, they hold whenever $D$ is in normal form and all local charges
 of $D$ vanish.
\end{enumerate}
\end{proposition}
\begin{proof}
For a pair of shears as above, \eqref{eq:shear-B} shows that $B$ is
invariant under $S\times S^\dual$, i.e. the pairing defined using $S$ is same as the pairing defined using base point, if and only if
\begin{equation}\label{eq:shear-compatible}
\begin{aligned}
B_{\mathrm{shear}}(u,v)-B_{\mathrm{base}}(u,v)
&=B(Sx,S^! y)-B(x,y)=-\big\langle\eta_D(\lambda^!),s(\lambda)\big\rangle
  -\big\langle\eta_{D^!}(\lambda),s^!(\lambda^!)\big\rangle=0.
\end{aligned}
\end{equation}
For all $\lambda, \lambda^!$.

If \eqref{eq:shear-compatible} holds then
$(SV_f)^{\perp_B}=S^\dual(V_f^{\perp_B})=S^\dual V^\dual_{f^\dual}$ by
Theorem~\ref{thm:C}, and the pairings transported by
$\Phi_{f_0}\circ S^{-1}$ and by $\Phi_{f_0}$ agree.  Conversely, if
$B(SV_f,S^\dual V^\dual_{f^\dual})=0$ for a single $f$, then, since
$B(V_f,V^\dual_{f^\dual})=0$ and $V_f$ projects onto $K_D$,
\eqref{eq:shear-B} forces \eqref{eq:shear-compatible}.  Also, the two
transported pairings agree only if $B$ is $S\times S^\dual$-invariant.

(1) Here $s=\frac1{|P_D|}\sum_{g}\delta^D_{f_0,g}$ and
$s^\dual=\frac1{|P_D|}\sum_{g}\delta^{D^\dual}_{f_0^\dual,g^\dual}$, we use $|P_D|=|P_{D^!}|$ since the
mirror bijection identifying the two fixed point sets, and
\eqref{eq:shear-compatible} is the average over $g\in P_D$ of
Theorem~\ref{thm:A} for the pairs $(f_0,g)$.

(2) Here $s=\wt(\,\cdot\,|_{f_0})^\circ$, $s^\dual=\wt(\,\cdot\,|_{f_0^\dual})^\circ$,
and \eqref{eq:shear-compatible} is \eqref{eq:GM-absolute} at $f_0$.  By
Theorem~\ref{thm:A}, \eqref{eq:GM-absolute} holds at one $f$ if and only if
it holds at all.  The proof of Corollary~\ref{cor:absolute} shows, for any
$D$ in normal form, that the left side of \eqref{eq:GM-absolute} is
$\frac12(\pi_z\otimes\pi_a)(\Omega_D)$ evaluated on $(\eta_D(\lambda^\dual),\lambda)$,
and that this vanishes when all local charges vanish.
\end{proof}

The two normalizations are not interchangeable in general.  For
$D=z_1z_2a_1a_2$ with D3 multiplicities $2,1,1,1$ on the intervals
$I_1,I_2,(a_1,a_2),(a_2,z_1)$. This a Nakajima quiver variety of type
$\widehat A_1$ with balanced D5 branes, but whose NS5 branes have
$\lch(z_1)=1$, $\lch(z_2)=-1$ respectively, direct computation shows \eqref{eq:Omega}
\[
 (\pi_z\otimes\pi_a)(\Omega_D)=(z_1-z_2)\otimes(a_1-a_2-b)\ne0.
\]
So the traceless subspaces $V^\circ_f$ and $V^{\dual\circ}_{f^\dual}$ are not
$B$-orthogonal, while the averaged ones always are.  

\begin{remark}[$SL_n$ versus $GL_n$]\label{rem:SL-GL} \rm
In \cite[\S\S7--8]{GM10} they consider $GL_n/P$ and Springer fibers, but the maximal torus is 
$H$ of $SL_n$ instead of the one of $GL_n$.  In general this is the traceless normalization, i.e. the framing torus
$\prod_i\C^*_{a_i}$ of $GL_n$ is replaced by its subtorus $\prod_ia_i=1$.  

For $G/B$ the
fixed points form a $W$-orbit, and the $W$-average of a weight
$\sum_ic_ia_i$ at fixed points of the torus of $GL_n$ is given by \[
\frac1{n!}\sum_{\tau\in S_n}\sum_{i=1}^n c_i a_{\tau(i)}
=\sum_{i=1}^n c_i\left(\frac1n\sum_{j=1}^n a_j\right)
=\frac{\sum_i c_i}{n}\sum_{j=1}^n a_j.
\]So
twisting by minus the average is the same as taking the traceless part.
Therefore the averaged and the traceless normalizations coincide.  This is the
computation ``$b\mapsto(b,b)$, $V=\bigcup_{\tau\in W}\{(\tau b,b)\}$'' of
\cite[\S7]{GM10}, whose average over $\tau$ is $b\mapsto(0,b)$. 

For general balanced
diagrams Proposition~\ref{prop:basepoint-free}(2) still gives GM duality
in the traceless normalization, with the same transported pairing.  For
charged diagrams it is either the {\em averaged normalization}, or a {\em base point choice}, that gives the correct pairing.
\end{remark}

\begin{remark}[Finite type $A$]\label{rem:typeA-normalization} \rm
For a circular diagram with an interval of multiplicity zero, i.e.\ a bow
variety of finite type $A$ \cite{rimanyi2020bow}, no tie winds around the
circle, and the parameters $b$ and $\kappa$ may be specialized to $0$.  With $b=\kappa=0$, \eqref{eq:Omega}
reads
\[
 \Omega_D=-2\,z_1\otimes\sum_i\lch(a_i)\,a_i+2\sum_j\lch(z_j)\,z_j\otimes a_1 ,
\]
which vanishes after the \emph{charge-weighted} normalization
$\sum_i\lch(a_i)a_i=0$, $\sum_j\lch(z_j)z_j=0$ (available when
$A=\sum_i\lch(a_i)\ne0$).  For $T^*(SL_n/P)$ presented with all D5 charges
equal, this is the traceless normalization, i.e.\ the maximal torus of
$SL_n$ of \cite[\S8]{GM10}.  
\end{remark}

\bigskip

These considerations leave us with the question: when is $\Phi_{f_0}$ an isomorphism?  Its restriction to $E_D$ is injective
because $r_f\circ\Phi_{f_0}$ is the identity on $E_D$.  The corresponding
map to ordinary cohomology is
\begin{equation}\label{eq:K-to-H2}
 K_D\to H^2(\X(D);\Q),\qquad \ell_j\mapsto c_1(\det\xi_{I_j}),
\end{equation}
Thus $\Phi_{f_0}$ is an isomorphism provided that:
\begin{enumerate}
\item[(a)] the determinant classes $c_1(\det\xi_{I_j})$, $1\le j\le m$, form
 a basis of $H^2(\X(D);\Q)$, and
\item[(b)] the sequence \label{sequence}$0\to E_D\to H^2_A(\X(D);\Q)\to H^2(\X(D);\Q)\to0$
 is exact.
\end{enumerate}
Under (b), the induced map $K_D\cong H_D/E_D\to H^2_A(\X(D);\Q)/E_D\to H^2(\X(D);\Q)$ is exactly
\eqref{eq:K-to-H2}, hence (a) is also necessary for $\Phi_{f_0}$ to be an
isomorphism.  Condition (b) holds, in particular, when $\X(D)$ is
equivariantly formal for $A$.  This is certainly not true for all bow varieties. It is shown by Y. Ji \cite{ji2023bow} injectivity of $H^*_A(\X(D),\Q)\to \oplus_{p\in\X(D)^A}H^*_A(p)$ can be failed.

Condition (a) is genuinely a condition.  It holds when $\X(D)$ is a Nakajima
quiver variety of affine type $A$ whose total $m$ tautological determinant classes
are linearly independent, by Kirwan surjectivity for quiver varieties
\cite{MN18}. It is worth mentioning that the same example by Y. Ji\cite{ji2023bow} shows that Kirwan surjectivity can fail for bow varieties.

Even if we assume Kirwan surjectivity, the independency reduces to a dimension count on $H^2$, which fails as well in some
degenerate cases.  For $\X(D)=\operatorname{Hilb}^1(\C^2)=\C^2$ one has
$m=1$ but $H^2=0$, so \eqref{eq:K-to-H2} is zero.  In this example the
tautological-model statement Theorem~\ref{thm:C} holds (trivially), while
Definition~\ref{def:GM} fails for the self-dual pair
$(\C^2,\C^2)$, since $K_p=0$ on both sides.  Although we do not pursue (a) further in this paper, we find it interesting:

\begin{question}
For which affine type $A$ brane diagrams is \eqref{eq:K-to-H2} an
isomorphism?  Equivalently, when is $H^2(\X(D);\Q)$ spanned by the
determinants of the tautological bundles, with no relations among the $m$
classes $\det\xi_{I_j}$?
\end{question}

\begin{remark}[Why the tautological model] \rm
From the point of view of 3d mirror symmetry the model $H_D=E_D\oplus K_D$
is the primary object rather than an approximation to $H^2_A$.  The elliptic stable
envelopes of \cite{BottaRimanyi} are sections of line bundles over an
extended elliptic cohomology scheme whose parameters are the equivariant
variables ($a_i$, $b$, $\h$) and the K\"ahler variables ($z_j$, $\kappa$),
one for each five-brane, regardless of the actual rank of $H^2(\X(D))$, and
the mirror parameter map $\varphi_D$ exchanges these parameter spaces. In this sense Theorem~\ref{thm:C} is GM duality in its natural
habitat, and we just consider Proposition~\ref{prop:GM-genuine} as the bridge to the
classical formulation.
\end{remark}

\begin{remark}\rm
Goresky--MacPherson's original examples are type $A$ bow mirror pairs,
namely $T^*(G/P)$ and its dual slice \cite{rimanyi2020bow}, and the Gale dual pair
$T^*\PP^{n-1}\leftrightarrow\widetilde{\C^2/\Z_n}$.
\end{remark}

\section{Weak order reversal}
\label{sec:weakOR}

\subsection{Chambers, attracting orders, weak order reversal}
\label{subsec:recollections-order}Fix a cocharacter $\sigma\colon\C^*\to\A_D$ such that $\X(D)^{\A_D}=\X(D)^{\sigma(\C^*)}$.  It defines the attracting sets
$\Att(f)=\{x\in\X(D):\lim_{t\to0}\sigma(t)x=f\}$ and the \emph{attracting
order} $\le$ on $P_D$, i.e. $g\in\overline{\Att(f)}\Rightarrow g\le f$, and we take the transitive closure. The order depends only on the
chamber of $\sigma$.

\begin{definition}[Weak order reversal]\label{def:weak-or}
The pair $(D,D^\dual)$, with chambers, has the \emph{weak order-reversal
property} if there is a partial order $\preceq$ on $P_D$ such that
$g\le h\Rightarrow g\preceq h$ (it refines the attracting order of $\X(D)$)
and $g^\dual\le h^\dual\Rightarrow h\preceq g$ (its transported opposite
refines the attracting order of $\X(D^\dual)$).
\end{definition}

\label{crit}The orders with weak order reversal property can be constructed from certain functions on the fixed points set,
$\Phi\colon P_D\to\R$ satisfies $f<g\Rightarrow\Phi(f)<\Phi(g)$ and
$f^\dual<g^\dual\Rightarrow\Phi(f)>\Phi(g)$, then any total order sorting
$\Phi$, is a partial order as in
Definition~\ref{def:weak-or}. And one can produce natural potential function using relative ample line bundle.

\begin{lemma}[Ample pairing]\label{lem:ample}
Let $L$ be an $\A_D$-linearized line bundle on $\X(D)$, relatively ample over
the affinization.  Then $h_{L,\sigma}(f)=\pair{\sigma}{\wt(L|_f)}$ is
strictly increasing along the attracting order, that is, $g<f\Rightarrow h(g)<h(f)$.
Changing the linearization changes $h$ by a constant.
\end{lemma}

\subsection{Positive determinant classes are relatively ample}
\begin{lemma}\label{lem:ampleness}
For every integral positive $\lambda\in K_D$, the line bundle $L_D(\lambda)$
is relatively ample over the affinization of $\X(D)$.
\end{lemma}
\begin{proof}
We work on the affine variety $\mathcal R_D$ of
Section~\ref{subsec:bow}, on which the conditions (S1), (S2) have
already been imposed. The quiver data defining each triangular part of $\mathcal{R}_D$ is given by:
\[
 \mathsf A\colon V_R\to V_L,\quad
 \alpha\colon V_R\to\C,\quad \beta\colon\C\to V_L,\qquad
 B_L\mathsf A-\mathsf A B_R+\beta\alpha=0,
\]
and denote $r=\dim V_R$, $l=\dim V_L$.
Here we forget the torus action with weight $t_1, t_2$. The map $\mathsf A$ has full rank by
\cite[Lemma~3.1]{nakajima2017cherkis}. Now we construct an invertible matrix from these data.

If $l>r$, denote $k=l-r$, then the square matrix,
\[
 M_+=\bigl[\mathsf A\ \ \beta\ \ B_L\beta\ \cdots\ \ B_L^{k-1}\beta\bigr]
\]
is invertible.  Starting with $W_0=\im\mathsf A$, successively adjoin
$\beta,B_L\beta,\dots$.  The moment-map equation implies
$B_L(\im\mathsf A)\subset\im\mathsf A+\im\beta$.
If this sequence of subspaces stops growing before reaching $V_L$, the
resulting proper subspace is $B_L$-invariant and contains
$\im\mathsf A+\im\beta$, contrary to (S2).  

If $r>l$, put $k=r-l$ and consider
\[
 M_-=\begin{bmatrix}
 \mathsf A\\ \alpha\\ \alpha B_R\\ \vdots\\ \alpha B_R^{k-1}
 \end{bmatrix}.
\]
This is invertible as well.  Indeed, intersect $\ker\mathsf A$ successively
with $\ker\alpha,\ker(\alpha B_R),\dots$.
This sequence will be stable if two consecutive intersections agree, if it is nonzero, then the moment-map
equation makes that intersection $B_R$ invariant and contained in
$\ker\mathsf A\cap\ker\alpha$, contradicting (S1).
Starting in dimension $k$, the intersection is therefore zero after $k$
steps.  

If $l=r$, use the invertible matrix $\mathsf A$ itself.

In either case its determinant is an invertible regular function $F_a$ on
the triangle locus, hence also on $\mathcal R_D$, with
\begin{equation}\label{eq:determinant-unit}
 F_a(gx)=\frac{\det(g_L)}{\det(g_R)}F_a(x).
\end{equation}
Recall the determinant line bundle $det\mathcal{V}_R$ and $det\mathcal{V}_L$ is defined by the trivial line bundle over $\mathcal{R}_D$
\[
L_R=\mathcal R_D\times\mathbb C,
\qquad
L_L=\mathcal R_D\times\mathbb C
\]
with gauge actions
\[
g\cdot_R(x,v)=(gx,\det(g_R)v),
\qquad
g\cdot_L(x,v)=(gx,\det(g_L)v).
\]
multiplication by $F_a$ defines an $G_D$-equivariant isomorphism
\[
\varphi:L_R\xrightarrow{\sim}L_L,
\qquad
\varphi(x,v)=(x,F_a(x)v).
\]
On the quotient this gives an isomorphism of the underlying line bundles.
Repeating this
argument across the D5 branes in an NS5 gap allows us to replace
the determinant character of any interval in that gap by that of
any other interval, without changing the isomorphism class of the
$G_D$-linearized line bundle.

\emph{The positive chamber.}
For positive integers $\lambda_j$, consider the character
\begin{equation}\label{eq:positive-GIT-character}
 \chi_\lambda(g)=\prod_{j=1}^m\det(g_{I_j})^{\lambda_j}.
\end{equation}
We are going to show all positive $\lambda$ share the same stable locus by showing that the same numerical criterion
\cite[Definition~2.7, Proposition~2.8]{nakajima2017cherkis} holds for any positive $\lambda$. We may replace $I_j$ by any representative interval in the same NS5 gap,
by \eqref{eq:determinant-unit}, so we can choose the rightmost one, then this is exactly the character of \cite[\S2.2(c)]{nakajima2017cherkis} with
$\vartheta_j=-\lambda_j$ (denoted $\nu_j^{\R}$ in that reference).

For an invariant graded subspace $S$ with $\alpha(S)=0$ and with
$\mathsf A\colon S_R\to S_L$ an isomorphism at every D5 brane, the condition is
\[
 -\sum_j\lambda_j\dim S_j\le0.
\]
Here $\dim S_j$ is constant along the $j$th NS5 gap because of the
isomorphism condition.  This inequality is automatic and is strict unless
$S=0$.  

For an invariant graded subspace $T$ containing all framing
images, with $\mathsf A\colon V_R/T_R\to V_L/T_L$ an isomorphism at every D5
brane, the condition is
\[
 -\sum_j\lambda_j\dim(V_j/T_j)\ge0.
\]
Since every $\lambda_j>0$, it holds exactly when $T=V$.
These conditions, and their strict versions, are independent of the
positive values of $\lambda_j$.  Hence all characters
\eqref{eq:positive-GIT-character} have the same semistable locus, which
equals their stable locus.

To identify this locus with the convention defining $\X(D)$, let $q_j$ be
the number of D3 intervals in the $j$th NS5 gap.  Repeating
\eqref{eq:determinant-unit} will identifies the linearization for
$\chi_0(g)=\prod_I\det(g_I)$ with that for
$\prod_j\det(g_{I_j})^{q_j}$, since $q_j$ is positive. Therefore we have,
\[
 \mathcal R_D/\!/_{\chi_\lambda}G_D
 =\mathcal R_D/\!/_{\chi_0}G_D=\X(D).
\]
Changing $\lambda$ only changes its polarization.

\emph{Ampleness.}
By affine GIT, the descended line bundle for $\chi_\lambda$ is relatively
ample for the projective morphism to
$Y_D=\operatorname{Spec}\C[\mathcal R_D]^{G_D}$.
With the convention of Section~\ref{subsec:bow}, this descended line is
exactly $L_D(\lambda)$.
\end{proof}

\subsection{The constant sum of heights}
Let $\lambda\in K_D$ and $\lambda^\dual\in K_{D^\dual}$ be integral and positive and
denote,
\begin{equation}\label{eq:sigmas}
 \sigma_D=-\eta_D(\lambda^\dual)\in E_D^*,
 \qquad
 \sigma_{D^\dual}=-\eta_{D^\dual}(\lambda)\in E_{D^\dual}^*,
\end{equation}
Choose positive integral multiples of $\lambda,\lambda^\dual$ if necessary so
that these rational cocharacters have integral lifts, and fix such lifts to
$\A_D$, $\A_{D^\dual}$.  Define the
\emph{heights}
\begin{equation}\label{eq:heights}
 H_D(f)=\bigl\langle\sigma_D,\wt(L_D(\lambda)|_f)\bigr\rangle,
 \qquad
 H_{D^\dual}(f^\dual)=\bigl\langle\sigma_{D^\dual},\wt(L_{D^\dual}(\lambda^\dual)|_{f^\dual})\bigr\rangle .
\end{equation}
The height $H_D(f)$ depends on the lift to $\A_D$, but it is easy to see the difference $H_D(f)-H_D(g)$ is lift-independent.
\begin{proposition}\label{prop:constant-sum}
For all $f,g\in P_D$,
\[
 H_D(f)-H_D(g)+H_{D^\dual}(f^\dual)-H_{D^\dual}(g^\dual)=0 ,
\]
so $H_D(f)+H_{D^\dual}(f^\dual)$ is independent of $f$.
\end{proposition}
\begin{proof}
$H_D(f)-H_D(g)=\pair{\sigma_D}{\wt(\lambda|_f)-\wt(\lambda|_g)}
=-\pair{\eta_D(\lambda^\dual)}{\delta^D_{f,g}(\lambda)}$. Similarly,
$H_{D^\dual}(f^\dual)-H_{D^\dual}(g^\dual)=-\pair{\eta_{D^\dual}(\lambda)}{\delta^{D^\dual}_{f^\dual,g^\dual}(\lambda^\dual)}$.
Then by relative GM identity \eqref{eq:GM} we conclude the proof.
\end{proof}

\subsection{The chamber computation}
In the rest of section $6$, we will verify the weak order reversal property for a particular chamber
\begin{equation}\label{eq:chamber}
 \sigma_b=-N\ll0,\qquad \sigma_{a_1}<\sigma_{a_2}<\dots<\sigma_{a_n},
\end{equation}
with $N$ much larger than all $|\sigma_{a_i}-\sigma_{a_{i'}}|$.
The mirror $\X(D^\dual)$ is equipped with the analogous chamber for mirror
D5 branes $a^\dual_1,\dots,a^\dual_m$ and the same cut.

Denote $d_i=R_i-R_{i+1}$ for $1\le i<n$ the number of NS5 branes between $a_i$
and $a_{i+1}$, and
\[
 d_\infty=m-\sum_{i<n}d_i=R_n+(m-R_1),
\]
the number of NS5 branes inbetween the $a_n$ and $a_1$ through the
cut.  Symmetrically, $d^\dual_j$ is the number of D5 branes between $z_j$ and
$z_{j+1}$, and $d^\dual_\infty=L_1+(n-L_m)$ is the number of D5 branes inbetween the $z_m$ and $z_1$ through the cut.

\begin{lemma}[Matched chambers]\label{lem:chambers}
Assume $d_\infty\ge1$.  Choose rational numbers $\varepsilon_i>0$ $(i<n)$ with
$\varepsilon=\sum_{i<n}\varepsilon_i<\frac1m$ (if $n=1$ take $\varepsilon=0$),
and $T>0$ such that
\begin{equation}\label{eq:theta-vector}
 \lambda^\dual_i=T\Bigl(\frac{d_i}m+\varepsilon_i\Bigr)\ \ (i<n),
 \qquad
 \lambda^\dual_n=T\Bigl(\frac{d_\infty}m-\varepsilon\Bigr)
\end{equation}
are integers and $T/m$ is integral, and denote $\lambda^\dual=\sum_i\lambda^\dual_i\ell^\dual_i$.  Then $\lambda^\dual$ is
positive, and $\sigma_D=-\eta_D(\lambda^\dual)$ satisfies
\begin{equation}\label{eq:sigma-props}
 \sigma_{t_1}=-\frac Tm<0,
 \qquad
 \sigma_{a_{i+1}}-\sigma_{a_i}=T\varepsilon_i>0,
 \qquad
 \sigma_{a_n}-\sigma_{a_1}=T\varepsilon<\frac Tm=|\sigma_{t_1}| .
\end{equation}
Also, $\sigma_D$ pairs nontrivially with every tangent weight at every
fixed point and lies in the chamber~\eqref{eq:chamber}. 
Moreover, the $\varepsilon_i$ can be chosen so that
$\X(D)^{\sigma_D(\C^*)}=\X(D)^{\A_D}$.  Symmetrically, if
$d^\dual_\infty\ge1$, then there is a positive $\lambda\in K_D$ such that
$\sigma_{D^\dual}=-\eta_{D^\dual}(\lambda)$ lies in the chamber
\eqref{eq:chamber} of $\X(D^\dual)$.
\end{lemma}
\begin{proof}
$\sum_i\lambda^\dual_i=T$, so by \eqref{eq:eta-t1} $\sigma_{t_1}=-T/m$ and
$\sigma_{a_{i+1}}-\sigma_{a_i}=\lambda^\dual_i-d_iT/m=T\varepsilon_i$.  Summing along all $i=1,...,n-1$,
$\sigma_{a_n}-\sigma_{a_1}=T\varepsilon<T/m$.

Since $\varepsilon<1/m\le d_\infty/m$, we have $\lambda^\dual_n=T(d_\infty/m-\varepsilon)>0$. Hence all $\lambda^!_i$ are positive.

From \eqref{eq:sigma-props} it is easy to see $\sigma_D$ are in the \eqref{eq:chamber}. We are going to show it pairs nontrivially with all tangent weight. Tangent weight at a fixed point is, modulo $\h$, $w=a_{i'}-a_i+ct_1$ with
$c\in\Z$, and $c\ne0$ if $i=i'$ (Section~\ref{subsec:E}).  By
\eqref{eq:tangent-pairing},
$\pair{\sigma_D}w=\sigma_{a_{i'}}-\sigma_{a_i}+c\,\sigma_{t_1}$.  

If $c\ne0$
then $|c\sigma_{t_1}|\ge T/m>|\sigma_{a_{i'}}-\sigma_{a_i}|$, hence $\langle \sigma_D,w\rangle\neq 0$. If $c=0$ then 
$\langle\sigma_D,w\rangle=(a_{i'}-a_i)\not=0$. So in either case the tangent weight is nonzero.

Since the chamber we choose is open, the standard argument shows that generic $\sigma_D$ in this chamber will satisfy $\X(D)^{\sigma_D(\C^*)}=\X(D)^{\A_D}$.

The mirror
statement is the same computation with the pairing coordinate given by the second half of
Lemma~\ref{lem:eta}.
\end{proof}

\begin{lemma}\label{lem:cut}
If the cut lies between two five-branes of opposite types, then
$d_\infty\ge1$ and $d^\dual_\infty\ge1$.
\end{lemma}
\begin{proof}
Follows from the definition of $d_{\infty}$ and $d_{\infty}^!$.
\end{proof}

\subsection{Weak order reversal}
\begin{theorem}[Weak order reversal]\label{thm:B}
Let $D$ be a circular brane diagram with $P_D\ne\varnothing$,
presented with a cut between two five-branes of opposite types, and equip
$\X(D)$ and $\X(D^\dual)$ with the chambers~\eqref{eq:chamber}.  Then
$(D,D^\dual)$ has the weak order-reversal property.

\end{theorem}
\begin{proof}
By Lemma~\ref{lem:cut}, $d_\infty\ge1$ and $d^\dual_\infty\ge1$, so
Lemma~\ref{lem:chambers} provides positive $\lambda^\dual\in K_{D^\dual}$ and
$\lambda\in K_D$ with $\sigma_D$, $\sigma_{D^\dual}$ in the two chambers.  By
Lemma~\ref{lem:ampleness}, $L_D(\lambda)$ and $L_{D^\dual}(\lambda^\dual)$ are
relatively ample, so by Lemma~\ref{lem:ample} $H_D$ is strictly increasing
along the attracting order of $\X(D)$ and $H_{D^\dual}$ along that of
$\X(D^\dual)$.  By Proposition~\ref{prop:constant-sum},
$H_D(f)=\mathrm{const}-H_{D^\dual}(f^\dual)$ is strictly decreasing along the
transported order.  Criterion~\ref{crit} applied to $\Phi=H_D$ finishes the
proof.
\end{proof}

\section{Relation to the Stable envelope mirror symmetry}
\label{sec:quadratic}
In this section we explain how to deduce relative Goresky-MacPherson duality from the existence of the mother function conjectured in \cite{AO21}. 
\subsection{Stable envelopes review}
Here we give a quick review of definition of elliptic stable envelope, for a more comprehensive discussion see \cite{BottaRimanyi}.

Given a bow diagram $D$, we consider a class 
\begin{multline*}
\alpha=\hbar\Bigg(
\sum_{i=1}^{n}\Big[
\operatorname{Hom}\!\left(\xi_{I_R(a_i)},\xi_{I_L(a_i)}\right)
+\operatorname{Hom}\!\left(\mathbb C_i,\xi_{I_L(a_i)}\right)
\Big]
+\sum_{j=1}^{m}t_2\,
\operatorname{Hom}\!\left(\xi_{I_L(z_j)},\xi_{I_R(z_j)}\right)
-\sum_{I}\operatorname{End}(\xi_I)
\Bigg)^\vee .
\end{multline*}
Here \(a_i\) runs over the D5 branes, \(z_j\) over the NS5 branes, and \(I\) over the D3 intervals.
This defines a line bundle $\Theta(\alpha)$ over $E_{\T_D}(X(D))$ via the elliptic Chern class $\Theta:\mathrm{K}_{\T_D}(\X(D))\to E_{\T_D}(\X(D)).$ 

Under the given cut with a chosen lift $\tilde{z}_j$ of K\"ahler parameters $z_j$, we define $U(\xi_j,\h^{c_j}\frac{\tilde{z}_j}{\tilde{z}_{j+1}} )$ for $j\leq m-1$, and $
U(\xi_m,\h^{c_m}\frac{\tilde{z}_m}{\tilde{z}_{1}\zeta} )$ when $j=m$, here $U(V,z)=\Theta\bigl((V-1^{\oplus rkV})(z-1)\bigr).$
Now we consider the line bundle
\[
\mathcal{L}_{\X(D)}:=\Theta(\alpha)\bigotimes_{j=1}^{m-1}U(\xi_{I_j},\h^{c_j}\frac{\tilde{z}_j}{\tilde{z}_{j+1}} )\otimes U(\xi_{I_m},\h^{c_m}\frac{\tilde{z}_m}{\tilde{z}_{1}\zeta} )
\]over $E_{\T_D\times \A_D^!}(\X(D))$
with
$c_j=-\lch(z_j)+N_j$, where $N_j$ is the number of D5 branes between the
consecutive lifted NS5 branes $\widetilde z_j,\widetilde z_{j+1}$, and with
$I_j$ the first D3 interval right of $z_j$ as in Section~\ref{subsec:K}. Let $e_f^-$ be the negative half of the tangent
weights at $f$ (with respect to the chamber \eqref{eq:chamber}).
\begin{definition}
The elliptic stable envelope associated to a chamber $\mathfrak{C}$ and fixed point $f$ in $\X(D)^{\T_D}$ is defined as a section $\mathrm{Stab}_D(f)$ of the bundle
\[
\bigl(\pi^*(i_f^*\mathcal{L}_{\X(D)}\otimes\Theta(-e^-_{f}))\bigr)^{-1}\boxtimes \mathcal{L}_{\X(D)}(\infty\Delta)
\] over $E_{\T_D\times \A_D^!}(\X(D))$, where $\pi:E_{\T_D\times \A_D^!}(\X(D))\to E_{\T_D\times \A_D^!}(f)$, such that it is supported on $\overline{Attr}^{\mathfrak{C}}_f$ and equal to $[Attr^{\mathfrak{C}}_f]$ when restrict to $f$.

Here $\Delta$ is the resonant locus inside $E_{\T_D/\A_D\times \A_D^!}(pt)$ defined in \cite[Sect~.2.3]{OkounkovInductive}, and  sections of $\mathcal{L}_{\X(D)}\bigl(\infty\Delta)$ can have arbitrary pole degree along $\Delta$.
\end{definition}
\subsection{Mother function and stable envelope mirror symmetry}
The stable envelope mirror symmetry is, up to sign,
\begin{equation}\label{eq:stable-envelope-mirror}
 \frac{\res{\Stab_D(f)}_g}{\res{\Stab_D(g)}_g}
 =\pm \varphi_D\left(
  \frac{\res{\Stab_{D^\dual}(g^\dual)}_{f^\dual}}
       {\res{\Stab_{D^\dual}(f^\dual)}_{f^\dual}}
 \right).
\end{equation}
Authors of \cite{AO21} propose a way to prove the stable envelope mirror symmetry by building up a mother function $m$ as a section of some line bundle $M_D$ over $E_{\T_D\times \A_D^!}(\X(D)\times \X(D^!))$ such that 
\[
m|_{g\times X(D^!)}=\Theta(e^-_{g})\varphi_D\bigl(\mathrm{Stab}_{D^!}(g^!)\bigr),\quad m|_{X(D)\times f^!}=\varphi_D\bigl(\Theta(e^-_{f^!})\bigr)\mathrm{Stab}_{D}(f).
\]
Therefore, the stable envelope mirror symmetry is simply the fact that the evaluation of $m$ at fixed point $(g, f^!)$ is order independent.

The existence of $M_D$ forces an isomorphism of line bundles
\[
\mathcal L_D(f)|_g\otimes\psi_D^*\Theta(e^-_{f^!})
\;\cong\;
\Theta(e^-_{g})\otimes\psi_D^*\!\left(\mathcal L_{D^!}(g^!)|_{f^!}\right)
\]
over $E_{\T_D\times \A_D^!}(g\times f^!)$. By taking the first Chern class we obtain an equality of quadratic forms
\[
Q_D(f)|_g-Q_D(g)|_g
=
\varphi_D\!\left(Q_{D^!}(g^!)|_{f^!}-Q_{D^!}(f^!)|_{f^!}\right).
\]
In the rest of this section we will show that when setting $\h=0$, such a quadratic form identity degenerates to the relative Goresky-MacPherson identity. 
\subsection{The Appell-Humbert quadratic forms for bow varieties}
The Appell-Humbert quadratic form determine the line bundle of which
$\Stab_D(f)$ is a section, which is defined
as follows. Recall that the line bundle of a stable envelope is defined by its $\alpha$ part and its $U$ part. The quadratic form associated to the $\alpha$ part is
\begin{align*}
 Q\alpha(a_i)&=\sum_{x_R\in X_{I_R}}\sum_{x_L\in X_{I_L}}(x_R-x_L+\h)^2
 +\sum_{x_L\in X_{I_L}}(a_i-x_L+\h)^2,\\
 Q\alpha(z_j)&=\sum_{x_L\in X_{I_L}}\sum_{x_R\in X_{I_R}}(x_L-x_R+t_1)^2,
 \qquad
 Q\alpha(I)=\sum_{x,x'\in X_I}(x-x'+\h)^2,
\end{align*}
and $Q\alpha_D=\sum_iQ\alpha(a_i)+\sum_jQ\alpha(z_j)-\sum_IQ\alpha(I)$.

The quadratic form associated to the $U$ part is given by
\begin{equation}\label{eq:QU}
 QU_D=2\sum_{j=1}^m(\beta_j+c_j\h)\sum_{x\in X_{I_j}}x .
\end{equation}
Put $\cW_D=Q\alpha_D+QU_D$.  Recall that $e_f^-$ is the negative half of the tangent
weights at $f$ (with respect to the chamber \eqref{eq:chamber}), and
$Q(e_f^-)=\sum_{w\in e_f^-}w^2$.  Then
\begin{equation}\label{eq:G}
 \cG_D(f)=\res{\cW_D}_f-Q(e_f^-),
 \qquad
 Q_D(f)=\cW_D-\cG_D(f),
 \qquad
 \Delta_D(f,g)=\cG_D(g)-\cG_D(f),
\end{equation}
and $Q_D(f)$ is the Appell--Humbert quadratic form of the line bundle of $\Stab_D(f)$.

The following is proved in the upcoming paper \cite{BottaRimanyi-Ahat}.
\begin{theorem}[Quadratic-form mirror identity \cite{BottaRimanyi-Ahat}]
\label{thm:BR-main}
For every circular brane diagram $D$ and all $f,g\in P_D$ we have
\begin{equation}\label{eq:3dMS-Q}
 \varphi_D\bigl(\Delta_{D^\dual}(f^\dual,g^\dual)\bigr)=-\Delta_D(f,g).
\end{equation}
\end{theorem}

\subsection{Mixed blocks}
In the coordinates $(a,b,\h,z,\kappa)$, a
quadratic form $\Theta$ in the parameters of $D$ has six blocks.  Its
\emph{mixed block} $\Mix_D\Theta$ is the part in the monomials $a_iz_j$,
$a_i\kappa$, $bz_j$, $b\kappa$. We regard $\Mix_D\Theta$ as an element of $Z_D\otimes E_D$, or
equivalently as a bilinear form on $E_D^*\times K_D$.

\begin{lemma}\label{lem:MixDelta}
For $f,g\in P_D$,
\begin{equation}\label{eq:MixDelta}
 \Mix_D\Delta_D(f,g)=-2\sum_{j=1}^m\beta_j\otimes\delta^D_{f,g}(\ell_j)
 =-\bigl(U_D(f)-U_D(g)\bigr)=-2\,\delta^D_{f,g}
 \quad\in Z_D\otimes E_D .
\end{equation}
\end{lemma}
\begin{proof}
The forms $Q\alpha_D$ and $Q(e_f^-)$ contain no K\"ahler variable, so only $QU_D$
contributes to the mixed block of $\cG_D(f)$, and the
terms $c_j\h$ in \eqref{eq:QU} do not.  Hence $\Mix_D(\cG_D(f))=2\sum_j\beta_j\otimes\ol{\wt(\ell_j|_f)}=U_D(f)$,
and $\Mix_D\Delta_D(f,g)=U_D(g)-U_D(f)$. Comparison with \eqref{eq:Wf} then completes the proof.
\end{proof}

\subsection{The other blocks}
\begin{proposition}\label{prop:blocks}
The $Z_D\otimes Z_D$
block and the $E_D\otimes E_D$ block of $\Delta_D(f,g)$ vanish. Consequently, the
quadratic-form mirror identity \eqref{eq:3dMS-Q} is equivalent to
Theorem~\ref{thm:A} together with the identity of the three blocks ($\h^2$, $\h(a\text{ or }b)$, $\h(z\text{ or }\kappa)$) involving
$\h$.
\end{proposition}
\begin{proof}
The $Z\otimes Z$ block vanishes because $QU_D$ is linear in the K\"ahler
variables and nothing else in $\cG_D(f)$ involves them.

For the $E\otimes E$ block we compare $Q\alpha_D|_f$ with the tangent
character at $\h=0$.  The tangent bundle of $\X(D)$ is
\begin{multline*}
 [T\X(D)]=\sum_i\bigl[(1-\h)\Hom(\xi_R,\xi_L)+\h\Hom(\xi_R,\C_i)+\Hom(\C_i,\xi_L)
 +\h\End(\xi_L)+\h\End(\xi_R)\bigr]_{a_i}\\
 +\sum_j\bigl[t_1\Hom(\xi_R,\xi_L)+t_2\Hom(\xi_L,\xi_R)\bigr]_{z_j}
 -(1+\h)\sum_I\End(\xi_I).
\end{multline*}
Here $t_1,t_2,\h$ denote multiplicative characters only here. Return now to additive weights, restrict to $\h=0$
(so the multiplicative character becomes $1$). Suppose
$T_fX=\bigoplus_w\mathbb C_w^{\oplus m_w}$ and consider the associated quadratic form $\sum_w m_w\,w^2.$  

At an NS5
brane the tangent weights are $\pm(x_L-x_R+\frac b2)$, contributing
$2Q\alpha(z_j)|_{\h=0}$.  At a D3 interval they are $x-x'$ twice, contributing $-2Q\alpha(I)|_{\h=0}$.  At a D5 brane the two copies of
$\Hom(\xi_R,\xi_L)$ cancel, and the remaining weights are $a_i-x_R$,
$x_L-a_i$, $x_L-x_L'$, $x_R-x_R'$.  

Denote $p_k$ for power sums, i.e. for $X=\{x_1,...,x_n\}$, $p_k(X)=\sum x_i^k$. Further denote $d=|X_R|$, $\ell=\lch(a_i)$, and using that at $\h=0$, $X_L-X_R$ consists of $\ell$ copies of $a_i$ (Lemma~\ref{lem:jump}), a direct
computation will show that both $Q\alpha(a_i)|_{\h=0}$ and half the
D5 tangent contribution equal
\[
 F=(2d+\ell+1)q-2s^2-2(\ell+1)a_i s+d(\ell+1)a_i^2.
\]
Here $s=p_1(X_R)$ and
$q=p_2(X_R)$.  Then $p_1(X_L)=s+\ell a_i$ and
$p_2(X_L)=q+\ell a_i^2$.  Since
the tangent weights come in pairs $\pm w$ on $\A_D$, we have $\sum_{w\in T_f\X(D)}\ol w^2=2Q\alpha_D|_{f,\h=0}$. Therefore $\bigl(Q\alpha_D-Q(e^-_f)\bigr)|_{E_D\times E_D}=0$, and $QU_D|_{E_D\times E_D}=0$ trivially hold. So the $E_D\times E_D$ block of $\Delta_D(f,g)$ vanishes.
\end{proof}

\section{Core--quotient combinatorics and 3d mirror symmetry}
\label{sec:content}

In this section we are going to show that the so-called core-quotient bijection for partitions is in fact $3d$ mirror symmetry of the torus fixed points of certain bow varieties. Moreover, we obtain a {\em colored} core-quotient formula via absolute GM identity.

\subsection{Charged partitions, their gluing}
The ordered pair $(\mu,e)$ where $\mu$ is a partition and $e\in \Z$ (the charge) is called a charged partition. Equivalently, a charged partition can be encoded by its B-set 
\[
B(\mu,e)=\{\mu_u-u+e :\  u=1,2,\ldots\} \subset \Z.
\]
The obtained set contains all $b\ll 0$ and no $b\gg 0$. Conversely, every subset of $\Z$  containing all $b\ll 0$ and no $b\gg 0$ is the B-set of a unique charged partition.

We define the conjugate of a charged partition $(\mu,e)$ to be $(\mu^t,-e)$. Its B-set is
\begin{equation}\label{eq:conj}
\{-1-b: \ b\not\in B(\mu,e)\}. 
\end{equation}

For a fixed $M$ and $M$ charged partitions $(\mu_1,e_1), (\mu_2,e_2), \ldots, (\mu_M,e_M)$ their $M$-gluing is defined to be the charged partition whose B-set is
\begin{equation*}
 \bigsqcup_{j=1}^M\{\,Mb+j-1:\ b\in B(\mu_j,e_j)\}
\qquad 
\left(
\text{in particular, } e=\sum_{j=1}^Me_j
\right).
\end{equation*}
If the glued charged partition is $(\RR,0)$, then $(\mu_1,…,\mu_M)$ is called the $M$-quotient of $\RR$.

Let us fix $M$. It is a fact that every charged partition $(\mu,e)$ is uniquely the $M$-gluing of $M$ charged partitions. 
Namely, if $B$ is the glued $B$-set, then $(\mu_j,e_j)$ has $B$-set
$\{b\in\Z: Mb+j-1\in B\}$ for $1\le j\le M$, see more details e.g. in \cite[Sect.~2.7]{JamesKerber}.

\begin{example} \rm
For $M=3$ the three charged partitions and their glued one 
\[
((3,1),2),\ 
((1),0),\ 
((2,2,1),1)\ 
\quad\leftrightarrow\quad
((10,7,5,4,3,2,1),3)
\]
determine each other via their B-sets.
\end{example}

The $M$-gluing of $(\emptyset,e_1),\ldots,(\emptyset,e_M)$ results in a charged partition from which an $M$-ribbon cannot be removed. Such partitions are called $M$-core partitions. 

If $(\mu,e)$ is the $M$-gluing of
$(\mu_1,e_1),\dots,(\mu_M,e_M)$, we write
$\operatorname{core}_M(\mu)$ for the underlying partition of the
$M$-gluing of $(\varnothing,e_1),\dots,(\varnothing,e_M)$.
This gluing has the same charge $e$ as $(\mu,e)$.

\subsection{Contents}
We can think of a partition $\mu$ as the Young-diagrams of its boxes $(u,v)$ with $1\leq v\leq \mu_u$. Define the {\em content} of a box $(u,v)$ of a charged partition $(\mu,e)$ to be $c(u,v)=v-u+e$.  When $M$ is fixed, a box with content congruent to $l$ mod $M$ is said to have color $l$. So, colors are indexed modulo $M$, i.e. color $M$ is color $0$.

For a multiset $X$ we define $N^M_l(X)$ to be the ``number of element of color $l$'', and $C^M_l(X)$ to be the {\em ``sum of the elements with color $l$''}. Evaluated on a charged partition, the boxes of the Young-diagram play the role of the multiset. For $((5,2),2)$ we have
\[
N^3_0=2, C^3_0=3+6=9, \qquad
N^3_1=2, C^3_1=4+1=5, \qquad
N^3_2=3, C^3_2=2+2+5=9.
\]

For $e=0$ and $M=1$ we recover the total content
$T(\mu)=C^1_0(\mu)$ of Section~\ref{subsec:Hilb}.

\subsection{Fixed points as arrays of charged partitions, their row and column gluing, 3dMS pairing}
\label{subsec:arrays}

Consider the diagram $D$ in normal form $z_1\cdots z_m\,a_1\cdots a_n$ with the cut between $a_n$ and $z_1$. Let $D^\dual$ be presented as in Lemma~\ref{lem:mirror}.  
We write $d=d_{I_m}$ for the multiplicity of the D3 interval between $z_m$ and $a_1$.
Following \cite{GyR} we will encode torus fixed points on $\X(D)$ by $n \times m$ arrays of charged partitions. 

Let $f\in P_D$ and fix $(i,j)$.  The function $p\mapsto\nu_f(p,i,j)$ has finite support.  It takes the values $0,1$ for $p\le 0$ and the values $0,-1$ for $p\ge1$.  Hence there is a unique charged partition $(\mu_{ij},e_{ij})=(\mu_{ij}(f),e_{ij}(f))$ such that,
\begin{equation}\label{eq:array}
 \nu_f(p,i,j)=\mathbf 1_{p\le0}-\mathbf 1_{p-1\in B(\mu_{ij},e_{ij})}
 \qquad(p\in\Z).
\end{equation}
We call $(\mu_{ij},e_{ij})_{1\le i\le n,\,1\le j\le m}$ the \emph{array of
$f$}.

The map $f\mapsto(\mu_{ij}(f),e_{ij}(f))$ is a bijection from $P_D$ to the
set of $n\times m$ arrays of charged partitions satisfying
\begin{equation}\label{eq:array-conditions}
 \sum_{j=1}^me_{ij}=-\lch(a_i)\ \ (1\le i\le n),\quad
 \sum_{i=1}^ne_{ij}=-\lch(z_j)\ \ (1\le j\le m),\quad
 \sum_{i,j}\Bigl(|\mu_{ij}|+\binom{e_{ij}}2\Bigr)=d,
\end{equation}
see details in \cite{GyR}. The sum $\sum
\binom{e_{ij}}{2}$ is the `core contribution' to $d$.
In the terminology of \cite[Sect.~5]{GyR}, the matrix $(e_{ij})$ is the
matrix of the core of the Maya diagram of $f$, the third condition in \eqref{eq:array-conditions} is the multiplicity of a core.

\begin{definition}\label{def:RK}
Let $f\in P_D$.  Let $(\RR_i(f),e_i)$ be the $m$-gluing of the $i$th row $(\mu_{i1},e_{i1}),\dots,(\mu_{im},e_{im})$ of
the array of $f$.  Let $(\KK_j(f),e^\dual_j)$ be the
$n$-gluing of the {\em conjugates} $(\mu_{1j}^t,-e_{1j}),\dots,(\mu_{nj}^t,-e_{nj})$ of the $j$th column.  
\end{definition}

When we evaluate contents on $\RR$ or $\KK$ below, we should always remember that they come equipped with the charges $e_i$ ($=-\lch(a_i)$) and $e^\dual_j$ ($=\lch(z_j)$), respectively.

Either of the lists, $\RR_i(f)$ or $\KK_j(f)$, determines the array, and hence $f$. The 3dMS pairing, given in terms of $\nu$ in \eqref{eq:N-mirror} translates to the array code of fixed points as follows, extending the transposition rule $\mu\mapsto\mu^t$ of Section~\ref{subsec:Hilb}:

\begin{theorem}[3d mirror symmetry is transposition of array]\label{thm:transpose}
For every $f\in P_D$, the array of the mirror fixed point $f^\dual\in
P_{D^\dual}$ is the conjugate transpose of the array of $f$, that is,
\[
 \mu_{ji}(f^\dual)=\mu_{ij}(f)^t,\qquad e_{ji}(f^\dual)=-e_{ij}(f)
 \qquad(1\le i\le n,\ 1\le j\le m).
\]
In particular, $\RR_j(f^\dual)=\KK_j(f)$ and
$\KK_i(f^\dual)=\RR_i(f)$. \qed
\end{theorem}

\subsection{Weights and contents}
\label{subsec:weights-contents}
Let us abbreviate $N_l=N^m_l$ and $C_l=C^m_l$ for the partitions $\RR_i$, and $N_k=N^n_k$ and $C_k=C^n_k$  for the partitions $\KK_j$. 
For $e\in\Z$ let $\mathrm{St}(e)$ define the ``staircase multiset''
\[
 \mathrm{St}(e)=\bigsqcup_{b=0}^{e-1}\{1,2,\dots,b\}\quad(e\ge0),
 \qquad
 \mathrm{St}(e)=\bigsqcup_{h=e}^{-1}\{h+1,\dots,-1,0\}\quad(e<0).
\]

\begin{lemma}\label{lem:weights}
Let $f\in P_D$ and $\RR_i=\RR_i(f)$.  For $1\le l\le m$,
\begin{equation}\label{eq:weights-contents}
 \ol{\wt(\ell_l|_f)}
 =\sum_{i=1}^n\Bigl[\bigl(N_l(\RR_i)+N_l(\mathrm{St}(e_i))\bigr)\,a_i
 +\frac b2\,\bigl(C_l(\RR_i)+C_l(\mathrm{St}(e_i))\bigr)\Bigr].
\end{equation}
\end{lemma}

In words, modulo $\h$ the weights of $\xi_{I_l}$ at $f$ are the characters
$a_i+\frac c2b$, one for each box of color $l$ and content $c$ in
$\RR_i$, together with some weights that depend only on the local
charges of the D5 branes.  If $\lch(a_i)=0$ for all $i$, there are no other
weights, and in particular $d_{I_l}=\sum_iN_l(\RR_i)$.

\begin{proof}
Let's use the following numbering of the lifted NS5 branes: number the lifted
NS5 branes along the universal cover so that $z_j^{(p)}$ has the number
$x=m(p-1)+j$.  The lifts with $x\le0$ are to the left of $a_i^{(0)}$, and
those with $x\ge1$ are to its right.  Then, using \eqref{eq:array}, there is a tie of $f$ between $a_i^{(0)}$ and the NS5 lift
with the number $x$ if and only if
\begin{equation}\label{eq:ties-R}
 x\le0\ \text{ and }\ x-1\notin B(\RR_i,e_i),
 \qquad\text{or}\qquad
 x\ge1\ \text{ and }\ x-1\in B(\RR_i,e_i).
\end{equation}
Let's also give the number $x$ to the lifted D3 interval
just to the right of the NS5 lift with the number $x$.  It is a lift of $I_l$
with $l\equiv x$ modulo $m$.  By Lemma~\ref{lem:intrinsic}, the portion of a
tie at $a_i^{(0)}$ over the interval with the number $x$ has weight
$a_i+\frac x2b$ modulo $\h$, both for $x\le0$ and for $x\ge1$.  A left tie
ending at the NS5 lift with the number $P\le0$ covers the intervals
$P,\dots,0$, and a right tie ending at the NS5 lift with the number $P\ge1$
covers the intervals $1,\dots,P-1$.  So, by \eqref{eq:ties-R} and by deck
transformations,
\[
 \ol{\wt(\ell_l|_f)}=\sum_{i=1}^n\
 \sum_{\substack{x\in X(B(\RR_i,e_i))\\ x\equiv l\ (\mathrm{mod}\ m)}}
 \Bigl(a_i+\frac x2\,b\Bigr),
 \qquad
 X(B)=\bigsqcup_{b\in B,\ b\ge0}\{1,\dots,b\}\ \sqcup
 \bigsqcup_{h\notin B,\ h<0}\{h+1,\dots,0\}.
\]
It remains to show that $X(B(\RR,e))$ is the union of the multiset of
the contents of $(\RR,e)$ and of $\mathrm{St}(e)$, for every charged
partition.  Let first $e=0$.  An element $b=\RR_u-u\ge0$ of
$B(\RR)$ contributes the contents of the boxes of the $u$th row that
lie to the right of the diagonal.  By\eqref{eq:conj}, an element $h<0$ not
in $B(\RR)$ has the form $h=-1-(\RR^t_v-v)$, and it contributes
the contents of the boxes of the $v$th column that lie on or below the
diagonal.  If $B$ has charge $e$, then $X(B+1)$ is obtained from the shifted
multiset $X(B)+1$ by adding $e$ copies of the element $1$ (removing $-e$
copies if $e<0$), and $\mathrm{St}(e+1)$ is obtained from $\mathrm{St}(e)+1$
in the same way.  So the claim follows for all $e$ by induction.
\end{proof}

\subsection{The content identity}
\label{subsec:content-identity}
 Observe that the multiplicity of the D3 interval of $D$ between $a_n$ and $z_1$,
which is the interval of $D^\dual$ between $z^\dual_n$ and $a^\dual_1$, is 
$d^\dual:=d-\sum_{i=1}^n\lch(a_i)$.

Our GM identity, rewritten in terms of an array of charged partitions, implies the following purely combinatorial theorem. 

\begin{theorem}[Content identity]\label{thm:content}
Let $f\in P_D$, and write $\RR_i=\RR_i(f)$,
$\KK_j=\KK_j(f)$.  Then
\begin{align}
 \sum_{j=1}^mC_k(\KK_j)+n\sum_{i\le k}N_0(\RR_i)
 +\frac nm\sum_{i=1}^nC_0(\RR_i)
 &=k\,d+\varkappa_k &&(1\le k\le n),\label{eq:content-K}\\
 \sum_{i=1}^nC_l(\RR_i)+m\sum_{j\le l}N_0(\KK_j)
 +\frac mn\sum_{j=1}^mC_0(\KK_j)
 &=l\,d^\dual+\varkappa^\dual_l &&(1\le l\le m),\label{eq:content-R}
\end{align}
where the numbers $\varkappa_k$ and $\varkappa^\dual_l$ depend only on the
local charges of $D$.  They are zero if all local charges of $D$ vanish.  In
particular the left sides of \eqref{eq:content-K} and \eqref{eq:content-R} do
not depend on $f\in P_D$.
\end{theorem}

\begin{proof}
It is enough to prove the first equality because the second one is the same for 
$D\to D^\dual$, $f\to f^\dual$. 

Consider the element
\[
 U_D(f)+\varphi_D\bigl(U_{D^\dual}(f^\dual)\bigr)\ \in\
 \operatorname{span}_\Q\{z_j,\kappa\}\otimes\operatorname{span}_\Q\{a_i,b\}
\]
of Section~\ref{subsec:absolute}.  By Theorem~\ref{thm:A} (in the form
\eqref{eq:GM-tensor}, see Section~\ref{subsec:proofA}) it does not depend on
$f\in P_D$. In Proposition~\ref{prop:absolute} we proved the common value to be $\Omega_D$.  Now we
evaluate both statements on the classes $\lambda=\ell_m$ and
$\lambda^\dual=\ell^\dual_k$.  Let $\varepsilon$ be the linear functional on
$\operatorname{span}_\Q\{z_j,\kappa\}$ with $\varepsilon(z_j)=0$ for all $j$
and $\varepsilon(\kappa)=-\frac1n$.  Let $\widetilde\sigma$ be the linear
functional on $\operatorname{span}_\Q\{a_i,b\}$ with $\widetilde\sigma(a_i)=1$
for $i\le k$, $\widetilde\sigma(a_i)=0$ for $i>k$, and
$\widetilde\sigma(b)=\frac2m$.  By \eqref{eq:beta}, $\varepsilon(\beta_l)$ is
$1$ for $l=m$ and $0$ for $l<m$.  So $\varepsilon$ restricts to the pairing
with $\ell_m$ on $Z_D$.  By \eqref{eq:u}, in which $A_i=a_i$ because $D$ is in
normal form, $\widetilde\sigma(\varphi_D(\beta^\dual_{k'}))$ is $1$ for $k'=k$
and $0$ otherwise.  By \eqref{eq:eta-char} this means that $\widetilde\sigma$
restricts to $\eta_D(\ell^\dual_k)$ on $E_D$.

Write $\widetilde N$ and $\widetilde C$ for the numbers $N$ and $C$ with the
terms with $\mathrm{St}$ of Lemma~\ref{lem:weights} added.  Recall that
$\varphi_D(a^\dual_j)=z_j$ and $\varphi_D(b^\dual)=-2\kappa$.  By
Lemma~\ref{lem:weights} for $D$ and for $D^\dual$,
\begin{align*}
 (\varepsilon\otimes\widetilde\sigma)\bigl(U_D(f)\bigr)
 &=2\,\widetilde\sigma\bigl(\ol{\wt(\ell_m|_f)}\bigr)
 =2\sum_{i\le k}\widetilde N_0(\RR_i)
 +\frac2m\sum_i\widetilde C_0(\RR_i),\\
 (\varepsilon\otimes\widetilde\sigma)\bigl(\varphi_D(U_{D^\dual}(f^\dual))\bigr)
 &=2\,\varepsilon\bigl(\varphi_D\bigl(\ol{\wt(\ell^\dual_k|_{f^\dual})}\bigr)\bigr)
 =\frac2n\sum_j\widetilde C_k(\KK_j).
\end{align*}
The sum of the two right sides is $\frac2n$ times the left side of
\eqref{eq:content-K}, with $N$ and $C$ replaced by $\widetilde N$ and
$\widetilde C$.  By Theorem~\ref{thm:A} it does not depend on $f$.  By
Proposition~\ref{prop:absolute} it is equal to
\[
 (\varepsilon\otimes\widetilde\sigma)(\Omega_D)
 =-\frac2n\sum_{i\le k}\lch(a_i)\bigl((i-n)+(k-i)\bigr)+\frac2n\,k\,d
 =\frac2n\Bigl(k\,d+(n-k)\sum_{i\le k}\lch(a_i)\Bigr).
\]
Here we used that $\varepsilon$ kills the $\varsigma_j$, and that
$\gamma(A-d)+cA=-2d\,\kappa\sum_ia_i-(A-d)\,b\sum_jz_j$. 

Finally, the corrections $\widetilde N-N$ and $\widetilde C-C$
come from the staircase multisets in Lemma~\ref{lem:weights}.
Their charges are $e_i=-\lch(a_i)$ (for $\RR_i$) and $e_j^\dual=\lch(z_j)$ (for $\KK_j$). Thus the corrections
depend only on the local charges of $D$ and vanish when all
local charges vanish. This proves
\eqref{eq:content-K}.
\end{proof}

\subsection{One D5 brane: the core--quotient formula}
\label{subsec:oneD5}

It is worth naming the content identity in the special case of $n=1$ D5 brane.  For simplicity denote
\[
 q=\lch(a),\qquad q_j=\lch(a^\dual_j)=-\lch(z_j)\quad(1\le j\le m),
 \qquad\text{so that}\qquad \sum_{j=1}^mq_j=-q .
\]
The array of $f\in P_D$ is a single row
$(\mu_1,q_1),\dots,(\mu_m,q_m)$, and
the number
\[
 \mathsf w=\sum_{j=1}^m|\mu_j|=d-\sum_{j=1}^m\binom{q_j}2
\]
is the same for all $f\in P_D$.  We write $\RR=\RR_1(f)$, it has charge $-q$.  The column partitions are
$\KK_j=\KK_j(f)=\mu_j^t$, with charges $-q_j$.  Let
$D-\mathsf w$ be the brane diagram obtained from $D$ by lowering all
multiplicities by $\mathsf w$ (with unchanged local charges).

\begin{corollary}[One D5 brane]\label{cor:oneD5}
Let $n=1$ and $f\in P_D$, with the notation above.
\begin{enumerate}
\item[(1)] The mirror fixed point $f^\dual$ is encoded by the $m$-tuple
 $(\KK_1,\dots,\KK_m)=(\mu_1^t,\dots,\mu_m^t)$.  The $m$-core of
 $\RR$ is the $m$-gluing of
 $(\varnothing,q_1),\dots,(\varnothing,q_m)$.  It depends only on $D$, and it
 is the row partition of the unique fixed point of $\X(D-\mathsf w)$.
 \item[(2)] For $1\le l\le m$, evaluate $C_l=C_l^m$ on both
$\RR$ and $\operatorname{core}_m(\RR)$ with charge $-q$,
so a box $(u,v)$ has content $v-u-q$ in either partition. Then
 \begin{equation}\label{eq:oneD5}
 \begin{aligned}
  C_l(\RR)-C_l(\operatorname{core}_m(\RR))
  &=-m\sum_{j=1}^m\bigl(T(\KK_j)-q_j|\KK_j|\bigr)
    +l\,\mathsf w-m\sum_{j\le l}|\KK_j|\\
  &=\phantom{-}m\sum_{j=1}^m\bigl(T(\mu_j)+q_j|\mu_j|\bigr)
    +l\,\mathsf w-m\sum_{j\le l}|\mu_j| .
 \end{aligned}
 \end{equation}
\item[(3)] Suppose that $q=0$.  Then $f\mapsto\RR_1(f)$ is a bijection
 from $P_D$ to the set of partitions with $d_{I_l}$ boxes of color $l$ for
 $1\le l\le m$.  Moreover $(\mu_1,\dots,\mu_m)$ is the $m$-quotient of
 $\RR$, we have $q_j=N_{j-1}(\RR)-N_j(\RR)$, and the
 following diagram commutes.
 \[
 \begin{tikzcd}[column sep=1.6cm,row sep=1.2cm]
  P_D \arrow[r,leftrightarrow,"\RR_1"]
      \arrow[d,leftrightarrow,"f\mapsto f^\dual"']
  & \bigl\{\text{partitions }\RR:\ N_l(\RR)=d_{I_l}\ (1\le l\le m)\bigr\}
      \arrow[d,leftrightarrow,"{\text{$m$-quotient, then transposition}}"]\\
  P_{D^\dual} \arrow[r,leftrightarrow,"(\RR_1{,}\dots{,}\RR_m)"']
  & \bigl\{(\KK_1,\dots,\KK_m):\
      \textstyle\sum_j|\KK_j|=d-\frac12\sum_jq_j^2\bigr\}
 \end{tikzcd}
 \]
 The `up-direction' of the right-hand arrow sends $(\KK_1,\dots,\KK_m)$
to the partition obtained by $m$-gluing
$(\KK_1^t,q_1),\dots,(\KK_m^t,q_m)$ in that order,
using the runner charges $q_j=-\lch(z_j)$ determined by $D$.
\end{enumerate}
\end{corollary}

So for one D5 brane, 3d mirror symmetry of fixed points is the core--quotient
bijection followed by the transposition of the components of the quotient.
The GM identity is the formula \eqref{eq:oneD5}, which expresses the content
sums of a partition, color by color, through its core and its quotient.

\begin{proof}
(1) The first claim is Theorem~\ref{thm:transpose} for $n=1$.  By
\eqref{eq:array-conditions}, the only array for $D-\mathsf w$ is
$(\varnothing,q_1),\dots,(\varnothing,q_m)$.

(2) Subtract \eqref{eq:content-R} for $D-\mathsf w$ and its fixed point from
\eqref{eq:content-R} for $D$ and $f$.  The local charges, and hence
$\varkappa^\dual_l$, are the same, and $d^\dual$ drops by $\mathsf w$.  For
$n=1$ we have $N_0(\KK_j)=|\KK_j|$ and
$C_0(\KK_j)=T(\KK_j)-q_j|\KK_j|$.  The second line follows
from $T(\mu^t)=-T(\mu)$.

(3) For $q=0$, Lemma~\ref{lem:weights} gives $N_l(\RR)=d_{I_l}$.
Conversely, the row of any partition $\RR$ is the array of a fixed
point of the diagram with one D5 brane and multiplicities $N_l(\RR)$.
The formula for $q_j$ is \eqref{eq:lch}, and the diagram commutes by (1).
\end{proof}

\begin{remark}[Sizes and contents]\label{rem:sizes} \rm
Let $n=1$ and $q=0$.  By Lemma~\ref{lem:weights},
$\ol{\wt(\ell_l|_f)}=N_l(\RR)\,a+\frac b2\,C_l(\RR)$, and the same
holds for the core.  Formula \eqref{eq:oneD5} compares the coefficients of
$b$ for $\RR$ and for its core.  The comparison of the coefficients of
$a$ is the comparison of the multiplicities of $D$ and of $D-\mathsf w$, that
is, $N_l(\RR)-N_l(\operatorname{core}_m(\RR))=\mathsf w$ for all
$l$.  Its sum over $l$ is the classical relation
$|\RR|-|\operatorname{core}_m(\RR)|=m\sum_j|\mu_j|$.  In this
sense \eqref{eq:oneD5} refines the classical relation between a partition,
its core and its quotient from sizes to contents.  If moreover all $q_j=0$,
then the core is empty, and \eqref{eq:oneD5} is the per-color content
identity for partitions with empty $m$-core.
\end{remark}

\begin{example}\label{ex:oneD5} \rm
Let $D=\ttt{$5$\fs $4$\fs $5$\fs $5$\bs}$.  Then $n=1$, $m=3$,
$(d_{I_1},d_{I_2},d_{I_3})=(4,5,5)$, $q=0$ and $(q_1,q_2,q_3)=(1,-1,0)$.  The
partition $\RR=(7,4,2,1)$ has $4$, $5$, $5$ boxes of colors $1$,
$2$, $3$, so it encodes a fixed point $f$ of $\X(D)$.  Its $3$-core is
$(1,1)$ and its $3$-quotient is
$(\mu_1,\mu_2,\mu_3)=((2),\varnothing,(1,1))$.  So $\mathsf w=4$, and the
mirror fixed point is encoded by
$(\KK_1,\KK_2,\KK_3)=((1,1),\varnothing,(2))$.  The content
sums $(C_1,C_2,C_3)$ are $(4,7,6)$ for $\RR$ and $(0,-1,0)$ for its
core.  On the right side of \eqref{eq:oneD5} the first term is
$-3\bigl((-1-2)+0+(1-0)\bigr)=6$, and the three cases $l=1,2,3$ of the
formula read
\[
 4-0=6+4-6,\qquad 7-(-1)=6+8-6,\qquad 6-0=6+12-12 .
\]
\end{example}

\bibliographystyle{alpha}
\bibliography{bow}

\begin{thebibliography}{BLPW16}

\bibitem[AO21]{AO21}
M.~Aganagic and A.~Okounkov.
\newblock Elliptic stable envelopes.
\newblock {\em Journal of the American Mathematical Society}, 34(1):79--133, 2021.

\bibitem[BD25]{BD25}
T.~Botta and H.~Dinkins.
\newblock Vertex functions for bow varieties and their mirror symmetry.
\newblock arxiv:2507.13219, 2025.

\bibitem[BFR23]{BFR24}
T.~Botta, A.~Foster, and R.~Rim{\'a}nyi.
\newblock Geometric {Bruhat} order on $(0,1)$-matrices, 2023.
\newblock arXiv:2311.05531.

\bibitem[BLPW16]{BLPW16}
T.~Braden, A.~Licata, N.~Proudfoot, and B.~Webster.
\newblock Quantizations of conical symplectic resolutions {II}: category {$\mathcal O$} and symplectic duality.
\newblock {\em Ast\'erisque}, 384:75--179, 2016.
\newblock With an appendix by I.~Losev.

\bibitem[BR25]{BottaRimanyi}
T.~M. Botta and R.~Rim\'anyi.
\newblock Bow varieties: {S}table envelopes and their 3d mirror symmetry.
\newblock {\em Duke Math.J.}, 174(13):2651--2748, 2025.

\bibitem[BR26]{BottaRimanyi-Ahat}
T.~M. Botta and R.~Rim\'anyi.
\newblock In preparation, 2026.

\bibitem[Che11]{cherkis2011instantons}
S.~Cherkis.
\newblock {Instantons on gravitons}.
\newblock {\em Communications in Mathematical Physics}, 306:449--483, 2011.

\bibitem[GM10]{GM10}
M.~Goresky and R.~MacPherson.
\newblock On the spectrum of the equivariant cohomology ring.
\newblock {\em Canad. J. Math.}, 62(2):262--283, 2010.

\bibitem[GR26]{GyR}
\'A. Gyenge and R.~Rim\'anyi.
\newblock Fixed point counts and motivic invariants of bow varieties of affine type {A}, with an {A}ppendix “{M}odularity” by {G}. {H}arcos.
\newblock to appear in Selecta Math., 2026.

\bibitem[Ji25]{ji2023bow}
Y.~Ji.
\newblock {Bow varieties as symplectic reductions of $T^*(GL_n/P)$}, with an appendix by {M}. {F}ranz.
\newblock {\em IMRN}, 2025(13):Paper No. rnaf194, 2025.

\bibitem[JK81]{JamesKerber}
G.~James and A.~Kerber.
\newblock {\em The Representation Theory of the Symmetric Group}, volume~16 of {\em Encyclopedia of Mathematics and its Applications}.
\newblock Addison--Wesley, 1981.

\bibitem[MN18]{MN18}
K.~McGerty and T.~Nevins.
\newblock Kirwan surjectivity for quiver varieties.
\newblock {\em Invent. Math.}, 212(1):161--187, 2018.

\bibitem[Nak23]{nakajima2018towards}
H.~Nakajima.
\newblock {Towards geometric Satake correspondence for Kac-Moody algebras--Cherkis bow varieties and affine Lie algebras of type $ A$}.
\newblock {\em Ann. Sci. Éc. Norm. Supér.}, 56(6):1777--1824, 2023.

\bibitem[NT17]{nakajima2017cherkis}
H.~Nakajima and Y.~Takayama.
\newblock {Cherkis bow varieties and Coulomb branches of quiver gauge theories of affine type {A}}.
\newblock {\em Selecta Mathematica}, 23:2553--2633, 2017.

\bibitem[Oko21]{OkounkovInductive}
A.~Okounkov.
\newblock Inductive construction of stable envelopes.
\newblock {\em Letters in Mathematical Physics}, 111(6):Paper No. 141, 2021.
\newblock 56 pp.

\bibitem[RS24]{rimanyi2020bow}
R.~Rim{\'a}nyi and Y.~Shou.
\newblock Bow varieties---geometry, combinatorics, characteristic classes.
\newblock {\em Comm. in Analysis and Geometry}, 32(2):507--575, 2024.

\bibitem[RSVZ22]{RSVZ22}
R.~Rim{\'a}nyi, A.~Smirnov, A.~Varchenko, and Z.~Zhou.
\newblock Three-dimensional mirror symmetry and elliptic stable envelopes.
\newblock {\em International Mathematics Research Notices}, 2022(13):10016--10094, 2022.

\bibitem[WY23]{WY23}
B.~Webster and P.~Yoo.
\newblock 3-dimensional mirror symmetry.
\newblock {\em Notices of the American Mathematical Society}, 70(9):1395--1406, 2023.

\end{thebibliography}

\end{document}